\documentclass[leqno]{article}
\usepackage{amssymb}
\usepackage{amsmath, enumerate, vmargin}
\usepackage[arrow, matrix]{xy}

\makeatletter
\renewcommand\section{\@startsection {section}{1}{\z@}%
 {-3.5ex \@plus -1ex \@minus -.2ex}%
 {2.3ex \@plus.2ex}%
 {\center \normalfont\large\bfseries}}
\makeatother

\setmarginsrb{3cm}{3cm}{3cm}{3cm}{1mm}{6mm}{0mm}{10mm}

\newtheorem{thm}{Theorem}[section]
\newtheorem{prop}[thm]{Proposition}
\newtheorem{cor}[thm]{Corollary}
\newtheorem{lem}[thm]{Lemma}
\newtheorem{defi}[thm]{Definition}
\newtheorem{remark}[thm]{Remark}
\newtheorem{example}[thm]{Example}
\newtheorem{pb}[thm]{Problem}

\newenvironment{rk}{\begin{remark}\rm}{\end{remark}}

\newenvironment{ex}{\begin{example}\rm}{\end{example}}

\newcommand{\real}{{\mathbb R}}
\newcommand{\nat}{{\mathbb N}}
\newcommand{\ent}{{\mathbb Z}}
\newcommand{\com}{{\mathbb C}}
\newcommand{\un}{1\mkern -4mu{\textrm l}}
\newcommand{\T}{{\mathbb T}}

\newcommand{\A}{{\mathcal A}}
\newcommand{\B}{{\mathcal B}}
\newcommand{\E}{{\mathcal E}}

\renewcommand{\H}{{\mathcal H}}

\newcommand{\M}{{\mathcal M}}
\newcommand{\N}{{\mathcal N}}
\newcommand{\R}{{\mathcal R}}
\renewcommand{\S}{{\mathcal S}}
\newcommand{\Z}{{\mathcal Z}}

\renewcommand{\a}{\alpha}
\renewcommand{\b}{\beta}
\newcommand{\g}{\gamma}
\renewcommand{\d}{\delta}

\newcommand{\e}{\varepsilon}
\renewcommand{\th}{\theta}
\renewcommand{\l}{\lambda}

\renewcommand{\e}{\varepsilon}
\newcommand{\f}{\varphi}
\renewcommand{\O}{{\Omega}}
\renewcommand{\o}{{\omega}}
\newcommand{\s}{\sigma}

\newcommand{\ma}{{\mathbb M}}
\newcommand{\Tr}{\mbox{\rm Tr}}
\newcommand{\tr}{\mbox{\rm tr}}
\newcommand{\ot}{\otimes}

\renewcommand{\t}{\tau}
\newcommand{\8}{\infty}
\newcommand{\el}{\ell}

\newcommand{\la}{\langle}
\newcommand{\ra}{\rangle}
\newcommand{\wt}{\widetilde}
\newcommand{\wh}{\widehat}
\newcommand{\n}{\noindent}
\newcommand{\pf}{\noindent{\it Proof.~~}}
\newcommand{\cqd}{\hfill$\Box$}
\newcommand{\be}{\begin{eqnarray*}}
\newcommand{\ee}{\end{eqnarray*}}
\newcommand{\beq}{\begin{equation}}
\newcommand{\eeq}{\end{equation}}

\font\ronde=cmmi8 at 14pt

 \newcommand\rM{{\mbox{\ronde m}}}
 \newcommand\rN{{\mbox{\ronde n}}}

\numberwithin{equation}{section}

\begin{document}


\title{A reduction method for noncommutative $L_p$-spaces
and applications }
\author{Uffe Haagerup\footnotemark\;\;
and Marius Junge\footnotemark\;\; and Quanhua Xu\footnotemark}

\date{}

\maketitle


\begin{abstract}

We consider the reduction of problems on general noncommutative
$L_p$-spaces to the corresponding ones on those associated with
finite von Neumann algebras. The main tool is an unpublished
result of the first-named author which approximates  any
noncommutative $L_p$-space by tracial ones.  We show that under
some natural conditions a map between two von Neumann algebras
extends to their crossed products by a locally compact abelian
group or to their noncommutative $L_p$-spaces. We present
applications of these results to the theory of noncommutative
martingale inequalities by reducing most recent general
noncommutative martingale/ergodic inequalities to those in the
tracial case.

\end{abstract}

 \setcounter{section}{-1}


 \makeatletter
 \renewcommand{\@makefntext}[1]{#1}
 \makeatother \footnotetext{\noindent
 $^*$ Partially supported by the Danish Natural Science
Research Council.\\
  $^\dagger$ Partially supported by the National Science
Foundation.\\
  $^\ddagger$ Partially supported by the Agence Nationale de Recherche.
 }



\section{Introduction}


The theory of noncommutative $L_p$-spaces has a long history going
back to pioneering works by von Neumann and Schatten
\cite{schatten-cr, schatten-id}, Dixmier \cite{dixLp}, Segal
\cite{segalLp} and Kunze \cite{kunze}. In the first constructions
the trace of a matrix or an operator replaces the integral of a
function, and the noncommutative $L_p$-spaces are composed of the
elements whose $p$-th power has finite trace. Later (around 1980),
generalizations to type III von Neumann algebras appeared due to
the efforts of Haagerup \cite{haag-Lp}, Hilsum \cite{hilsum},
Araki and Masuda \cite{am}, Kosaki \cite{kos-int} and Terp
\cite{terp-int}. These algebras have no trace, and therefore the
integration theory has to be entirely redone. These
generalizations were motivated and made possible by the great
progress in operator algebra theory, in particular the
Tomita-Takesaki theory and Connes's spectacular results on the
classification of type III factors.

Since the early nineties and the arrival of new theories such as
those of operator spaces and free probability, noncommutative
integration has been living another period of stimulating new
developments. In particular, noncommutative Khintchine and
martingale inequalities have opened new perspectives. It is well
known nowadays that the theory of noncommutative $L_p$-spaces is
intimately related with many other fields such as Banach spaces,
operator algebras, operator spaces, quantum probability and
noncommutative harmonic analysis. Although they correspond to
separate research directions, these fields  present many links and
interactions through  noncommutative $L_p$-spaces. For instance,
Pisier's works, in particular his theory  on vector-valued
noncommutative $L_p$-spaces \cite{pis-ast}, have started to
exhibit these interactions. Since the establishment of the
noncommutative Burkholder-Gundy inequalities in \cite{px-BG}, many
classical martingale inequalities have been transferred to the
noncommutative setting. These include the Doob inequality in
\cite{ju-doob}, the Burkholder/Rosenthal inequalities in
\cite{jx-burk, jx-ros}, the weak type $(1, 1)$ inequality for
martingale transforms in \cite{ran-mtrans} and the Gundy
decomposition in \cite{par-ran-gu}. This rapid development of
noncommutative inequalities is largely motivated by quantum
probability and operator space theory. Notably, the latter theory
has inspired many new ideas and provides numerous tools. These
noncommutative inequalities have, in return, interesting
applications to operator space theory, quantum probability and
more generally noncommutative analysis.

The recent works on the complete embedding of Pisier's operator
Hilbert space $OH$ into a noncommutative $L_1$ in \cite{ju-OH} and
certain Hilbertian homogeneous operator spaces into noncommutative
$L_p$-spaces in \cite{xu-embed} show well how useful these
$L_p$-spaces are as tools and examples for operator spaces. These
works also illustrate the power of quantum probabilistic methods
in operator spaces since Khintchine-type inequalities for free
random variables play a key role there. This illustration is
further shown in the works \cite{pisshlyak} and \cite{xu-gro} on
operator space Grothendieck theorems. It should be emphasized that
type III von Neumann algebras are unavoidable in all these works,
for Pisier \cite{pis-III} showed that $OH$ cannot completely embed
in a semifinite $L_1$.

While the tracial noncommutative $L_p$-spaces are a rather
transparent generalization of the usual $L_p$-spaces, all existing
(equivalent) constructions of type III $L_p$-spaces are quite
heavy and based on the Tomita-Takesaki theory. This makes them
much less pleasant and handy. For instance, the  $L_p$-spaces
associated to a type III algebra do not form a real interpolation
scale (cf. \cite{px-survey}). Another main drawback is the lack of
a reasonable analogue of weak $L_1$ in the type III case. It is
however well known that weak $L_1$-spaces are of paramount
importance in analysis notably through the Marcinkiewicz
interpolation theorem.

\smallskip

It is thus desirable to reduce or approximate type III
$L_p$-spaces to or by semifinite ones. This is exactly the
objective of a unpublished work \cite{haag-red} of the first named
author about three decades ago. The result there can be stated as
follows. Given a $\s$-finite von Neumann algebra $\M$ equipped
with a normal faithful state $\f$ there exists another $\s$-finite
von Neumann algebra $\R$ and a normal faithful state $\wh\f$ on
$\R$ verifying the following properties:
 \begin{enumerate}[(i)]

\item $\M$ is a von Neumann subalgebra of $\R$, the restriction of
$\wh\f$ to $\M$ is equal to $\f$ and there exists a
state-preserving normal faithful conditional expectation from $\R$
onto $\M$;

\item there exists an increasing sequence $(\R_n)_n$ of finite von
Neumann subalgebras of $\R$ such that their union is w*-dense in
$\R$ and such that each $\R_n$ is the image of a state-preserving
normal faithful conditional expectation.
 \end{enumerate}
Property (i) allows us to view $L_p(\M)$ isometrically as a
subspace of $L_p(\R)$. (ii) insures  that the sequence
$(L_p(\R_n))_n$ is increasing and $\bigcup_nL_p(\R_n)$ is dense in
$L_p(\R)$ for $p<\8$. This implies that $L_p(\R)$, so $L_p(\M)$,
can be approximated by the $L_p(\R_n)$, which are based on finite
algebras. This approximation theorem reduces many geometrical
properties of general noncommutative $L_p$-spaces to the
corresponding ones in the tracial case. This is indeed the case
for  all those properties which are of a finite-dimensional nature
such as Clarkson's inequalities, uniform convexity/smoothness and
type/cotype.

The preceding reduction theorem  has found more and more
applications since the new developments of noncommutative
$L_p$-spaces mentioned previously. In fact, it plays a crucial
role in many recent works. See, for instance, \cite{jios},
\cite{jupar-ros}, \cite{jupar-mix}, \cite{jx-erg}, \cite{xu-max}
and \cite{xu-descrip}. Moreover, in these works one needs the
precise form of $\R$ and $\R_n$ as constructed in \cite{haag-red}.
On the other hand,  applications of this theorem often require
additional results such as those on the extension of maps on von
Neumann algebras to their crossed products or noncommutative
$L_p$-spaces.

\smallskip

The manuscript \cite{haag-red}, however, has been circulated only
in a very limited circle of people. We feel that it would be
helpful to make it accessible to the general public. This explains
why we decide to include the proof of the previous reduction
theorem following the presentation of \cite{haag-red}. This
reproduction corresponds to section~\ref{Reduction} below. In the
second part of this article we show well how useful this theorem
is for noncommutative martingale and ergodic inequalities. The
first part contains two extension results for maps on von Neumann
algebras, one to their crossed products and another to their
noncommutative $L_p$-spaces. These extension results are of
interest for their own right.

\smallskip

The paper is organized as follows. In section~1 we summarize
necessary preliminaries on crossed products and noncommutative
$L_p$-spaces.   For these $L_p$-spaces we use the construction
\cite{haag-Lp} of the first-named author. Today, they are commonly
called Haagerup noncommutative $L_p$-spaces. In
section~\ref{Reduction} we prove the reduction theorem mentioned
previously. Our presentation follows the unpublished manuscript
\cite{haag-red}. The tool is crossed products.
Section~\ref{Applications to noncommutative Lp-spaces} presents
the first application of the reduction theorem to noncommutative
$L_p$-spaces. The result there, already quoted previously, says
roughly that a type III noncommutative $L_p$-space with $p<\8$ can
be approximated by tracial ones. This is in fact the original
intention of \cite{haag-red}. In section~\ref{Extensions of maps
to crossed products} we deal with the extension of a map between
two von Neumann algebras to their crossed products by a locally
compact abelian group. The extension of such a map to the
corresponding $L_p$-spaces is treated  in section~\ref{Extensions
of maps to noncommutative Lp-spaces}. The second extension is much
subtler than the first one. Its proof involves Kosaki-Terp's
interpolation theorem. Note, however, that this extension to $L_p$
is quite obvious in the tracial case. The last two sections
contain applications to martingale/ergodic inequalities. The first
one is devoted to square function type inequalities and the second
to maximal inequalities.  Most results in these two sections are
not new and some arguments also exist in the literature. For some
results there the reduction to the tracial case is not really
necessary. But for some others we do not know other methods than
the reduction, for example, for the maximal ergodic inequalities.
We feel that it would be helpful for the reader to have a complete
picture of how to reduce these inequalities to the tracial case.
This is why all inequalities in consideration are properly stated
and some known arguments are also included.


\section{Preliminaries}
 \label{Preliminaries}


In this section we collect some necessary preliminaries on crossed
products and noncommutative $L_p$-spaces used throughout the
paper.


\subsection{Crossed products}
 \label{Crossed products}


Our references for modular theory and crossed products are
\cite{kar-II}, \cite{ped-tak}, \cite{stratila}, \cite{tak-II} and
\cite{tak-cross}. Throughout this paper $\M$ will always denote a
von Neumann algebra acting on a Hilbert space $H$. Let $G$ be a
locally compact abelian group equipped with Haar measure $dg$, and
$\wh G$ its dual group equipped with Haar measure $d\g$. We choose
$dg$ and $d\g$ so that the Fourier inversion theorem holds. Let
$\alpha$ be a continuous automorphic representation of $G$ on
$\M$. We denote by $\R=\M\rtimes_\a G$ the crossed product of $\M$
by $G$ with respect to $\a$. Recall that $\R$ is the von Neumann
algebra on $L_2(G, H)$ generated by the operators $\pi_\a(x)$,
$x\in \M$ and $\l(g)$, $g\in G$, which are defined by
 $$\big(\pi_{\a} (x) \xi \big) (h)=\a^{-1}_{h}(x)
 \xi (h),\quad
 \big(\l(g) \xi\big) (h) = \xi (h-g), \quad  \xi\in L_2(G, H),\;
 h\in G.$$
These operators  satisfy the following commutation relation:
 \beq\label{com pi-l}
 \l(g) \pi(x) \l(g)^{\ast}=\pi(\a_{g}(x)),\quad
 x\in\M, \; g\in G.
 \eeq
Consequently, the family of all linear combinations of
$\pi_\a(x)\l(g)$, $x\in\M$, $g\in G$, is a w*-dense involutive
subalgebra of $\R$. Recall that $\pi_{\alpha}$ is a normal
faithful representation of $\M$ on $L_2(G, H)$, so we can identify
$\M$ and $\pi_{\alpha}(\M)$. In the sequel, we will drop the
subscript $\a$ from $\pi_\a$ whenever no confusion can occur.

The action $\a$ admits a dual action $\wh\a$ of the dual group
$\widehat G$ on the crossed product $\R$ defined as follows. Let
$w$ be  the following unitary representation of $\wh G$ on $L_2(G,
H)$:
 $$\big(w(\g)\xi\big)(h) = \overline{\g(h)}\xi(h),\quad
 \xi\in L_2(G, H),\; h \in G, \ \g \in\wh G.$$
Then the dual action $\wh\a$ is implemented by $w$:
 $$\wh\a_{\g} (x) = w(\g)xw(\g)^{\ast},
 \quad x \in \R, \ \g \in \wh G.$$
It is easy to check that
 \beq\label{dual action}
 \wh\a_{\g}(\pi(x)) = \pi(x),\quad
 \wh\a_{\g}(\l(g)) = \overline{\g(g)}\l(g),
 \quad x\in\M,\ g \in G, \ \g\in\wh G.
 \eeq
Recall that $\pi(\M)$ is the algebra of the fixed points of
$\wh\a$. Namely,
 \beq\label{invariance of M by dual action}
 \pi(\M)=\{ x\in \R\;:\;
 \wh\a_{\g}(x) = x,\ \forall\; \g\in \wh G\}.
 \eeq
There exists a normal semifinite faithful ({\it n.s.f.} for short)
operator-valued weight $\Phi$ in the sense of \cite{haag-ovwI}
from $\R$ to $\pi(\M)$ defined by
 \beq\label{o-v weight}
 \Phi(x)=\int_{\wh G}\wh \a_{\g} (x)\,d\g,\quad x\in
 \R_+.
 \eeq
Moreover,  $\Phi$ is $\widehat\alpha$-invariant, i.e.,
 $$\Phi\circ\wh\a_{\g}=\Phi,\quad \g\in\wh G.$$
Note that $\Phi$ takes values in the extended positive part of
$\pi(\M)$ and can be defined on the extended positive part of $\R$
too.  We can easily determine the action of $\Phi$ on the elements
in a w*-dense involutive subalgebra of $\R$. Indeed, let $f\,:\,
G\to\M$ be a compactly supported w*-continuous function. Then
 $$x_f=\int_G\pi(f(g))\l(g)\,dg$$
defines an operator in $\R$. One can check that the family of all
such operators $x_f$ forms a w*-dense involutive subalgebra of
$\R$. If additionally $x_f\ge 0$, then
 \beq\label{o-v weight bis}
 \Phi(x_f)=\pi(f(0)).
 \eeq
See \cite{haag-ovwI} for more details.

Let $\f$ be a normal semifinite weight on $\M$. Then $\f$ admits a
dual weight $\wh\f$ on the crossed product $\R$ given by
  \beq\label{dual weight}
  \wh\f= \f\circ\pi^{-1}\circ \Phi.
  \eeq
By \eqref{o-v weight bis} for every positive $x_f$ as above we
have
 $$\wh\f(x_f)=\f(f(0)).$$
$\wh\f$ is normal and semifinite; if $\f$ is faithful, so is
$\wh\f$. Since $\Phi$ is $\wh\a$-invariant, by \eqref{dual
weight}, $\wh\f$ is $\wh\a$-invariant too:
 $$\wh\f\circ\wh\a_{\g}=\wh\f,\quad\g\in\wh G.$$
Moreover, the map $\f\mapsto\wh\f$ is a bijection from the set of
all normal semifinite weights on $\M$ onto the set of all normal
semifinite $\wh\a$-invariant weights on $\R$ (cf.
\cite[section~19.8]{stratila}). The modular automorphism group of
the dual weight $\wh\f$ is uniquely determined by
 \beq\label{modular gp of dual weight}
 \s_{t}^{\wh\f} (\pi(x)) = \pi(\s_{t}^{\f}(x)),\quad
 \s_{t}^{\wh\f}(\l(g)) = \l(g),\quad x \in\M,\
 g \in G,\  t \in\real.
 \eeq
Thus $\s_{t}^{\wh\f}$ leaves all $\l(g)$ invariant, and if we
identify $\pi(\M)$ with $\M$, the restriction of $\s_{t}^{\wh\f}$
on $\M$ coincides with $\s^{\f}_{t}$.

\smallskip

The case where $G$ is a discrete group will play an important role
later. In this case the previous construction becomes simpler. In
particular, the operator-valued weight $\Phi$ defined by
\eqref{o-v weight} becomes a conditional expectation, which is
uniquely determined by
  \begin{eqnarray}\label{cond expect}
  \Phi \big(\l(g)\pi(x)\big) = \left\{\begin{array}{ll}
  \pi(x) & \textrm{ if } g=0, \\
  0  &\textrm{ if } g\neq 0,
  \end{array}\right. \quad x\in\M,\; g\in G.
  \end{eqnarray}
On the other hand, if $\f$ is a normal (faithful) state on $\M$,
$\wh\f$ is a normal (faithful) state on $\R$ determined by
 \begin{eqnarray}\label{dual state}
  \wh\f\big(\l(g)\pi(x)\big) = \left\{\begin{array}{ll}
  \f(x) & \textrm{ if } g=0, \\
  0  &\textrm{ if } g\neq 0,
  \end{array}\right.  \quad x\in\M,\; g\in G.
  \end{eqnarray}


\subsection{Noncommutative $L_p$-spaces}


We now introduce Haagerup noncommutative $L_p$-spaces following
\cite{haag-Lp} and \cite{terp}. All results quoted below without
extra reference can be found in these two papers. Throughout this
subsection, $\f$ will be a fixed {\it n.s.f.} weight on $\M$, and
$\R$ will denote the crossed product $\M\rtimes_{\s}\real$ of $\M$
by $\real$ with respect to the modular automorphism group
$\s=\s^\f$. We will identify $\M$ and the subalgebra $\pi(\M)$ of
$\R$. It is well known that $\R$ is semifinite  and there exists a
unique {\it n.s.f.} trace $\t$ on $\R$ such that
 \beq\label{unique trace}
 (D\wh\f\;:\;D\t)_t=\l(t), \quad t\in\real,
 \eeq
where $(D\wh\f\;:\;D\t)_t$ denotes the Radon-Nikod\'ym cocycle of
$\wh\f$ with respect to $\t$.  We will call this trace the
canonical {\it n.s.f.} trace on $\R$. Moreover, $\t$ is the unique
{\it n.s.f.} trace on $\R$ satisfying
 \beq\label{relation tau-sigma}
 \t\circ\wh\s_{t} = e^{-t}\t,\quad t\in\real.
 \eeq
Given a normal semifinite weight $\psi$ on $\M$, let $h_\psi$
denote the Radon-Nikod\'ym derivative of the dual weight $\wh\psi$
with respect to $\t$, which is the unique positive selfadjoint
operator affiliated with $\R$ such that
 $$\wh\psi(x)=\t(h_\psi x),\quad x\in \R_+.$$
Here $\t(h_\psi x)$ is understood as
$\t(h_\psi^{1/2}xh_\psi^{1/2})$.  Then by \eqref{relation
tau-sigma}, we have
 \beq\label{L1 under dual action}
 \wh\s_{t}(h_\psi) = e^{-t}h_\psi,\quad t\in\real.
 \eeq
Recall that the map $\psi\mapsto h_\psi$ is a bijection from the
set of all normal semifinite weights on $\M$ onto the set of all
positive selfadjoint operators affiliated with $\R$ satisfying
\eqref{L1 under dual action} (cf. \cite[Proposition~II.4]{terp}).

In particular, the dual weight $\wh\f$ of our distinguished weight
$\f$ has a Radon-Nikod\'ym derivative $D_\f$ with respect to $\t$.
We will call $D_\f$ the density operator of $\f$ and will often
denote it by $D$ whenever no confusion can occur. Then by
\eqref{unique trace} the regular representation $\l(t)$ of $\real$
on $L_2(\real, H)$ is given by
 $$\l(t)=D^{it},\quad t\in\real.$$

\smallskip

Now we are ready to define  Haagerup noncommutative $L_p$-spaces.
Let $L_0(\R, \t)$ denote the topological involutive algebra of all
operators on $L_2(\real, H)$ measurable with respect to $(\R, \t)$
(cf. \cite{nelson} and \cite[Chapter I]{terp}). Let
$0<p\leq\infty$. Then the Haagerup $L_p$-space associated with
$(\M,\f)$ is defined to be
 $$L_p(\M, \f)=\big\{x\in L_0(\R, \t)\, :\,
 \wh\s_{t}(x)=e^{-t/p}x,\ \forall\; t \in \real \big\}.$$
The spaces $L_p(\M,\f)$  are closed selfadjoint linear subspaces
of $L_0(\R,\tau)$. It is not hard to show that $L_\infty(\M,\f)$
coincides with $\M$. On the other hand, it is well known that the
map $\o\mapsto h_\o$ on $\M_*^+$ extends to  a linear
homeomorphism from $\M_*$ onto $L_1(\M,\f)$ (equipped with the
vector space topology inherited from $L_0(\R,\tau) $). This
permits to transfer the norm of $\M_*$ into a norm on
$L_1(\M,\f)$, denoted by $\|\;\|_1$. Moreover, $L_1(\M,\f)$ is
equipped with a distinguished contractive positive linear
functional $\tr$, the ``trace'', defined by
 $$\tr\,(h_\o)=\o(1),\quad \o\in \M_*.$$
Consequently, $\|h\|_1=\tr\,(|h|)$ for every $h\in L_1(\M,\f)$.

\smallskip

Let $0<p<\infty$. If $x=u|x|$ is the polar decomposition of $x\in
L_0(\R,\tau)$, then
 $$x\in L_p(\M,\f)\Leftrightarrow u\in \M
 \mbox{ and }|x|\in L_p(\M, \f)\Leftrightarrow u\in \M
 \mbox{ and }|x|^p\in L_1(\M, \f).$$
For $x\in L_p(\M,\f)$ set $\|x\|_p=\|\,|x|^p\|_1^{1/p}$. Then
$\|\;\|_p$ is a norm or a $p$-norm according to  $1\leq p<\infty$,
or $0<p<1$. The associated vector space topology coincides with
that inherited from $L_0(\R,\tau)$.

Another important link between the spaces $L_p(\M,\f)$ is the
external product: in fact, the product of $L_0(\R,\tau)$,
$(x,y)\mapsto xy$, restricts to a bounded bilinear map
$L_p(\M,\f)\times L_q(\M,\f)\to L_r(\M,\f)$, where $1/r=1/p+1/q$.
This bilinear map has norm one, which amounts to saying that the
usual H\"older inequality extends to Haagerup $L_p$-spaces. In
particular, if $1/p+1/p'=1$, then the bilinear form $(x,y)\mapsto
\tr\,(xy)$ defines a duality bracket between $L_p(\M)$ and
$L_{p'}(\M,\f)$, for which $L_{p'}(\M,\f)$ coincides
(isometrically) with the dual of $L_p(\M,\f)$ (if $p\neq\infty$);
moreover we have the tracial property:
 $$\tr(xy)=\tr\,(yx),\quad x\in L_p(\M,\f),\; y\in L_{p'}(\M,\f).$$

\smallskip

The distinguished weight $\f$ can be recovered from the tracial
functional $\tr$. Let
 $$\rN_\f=\{x\in\M\;:\; \f(x^*x)<\8\},\quad
 \rM_\f=\rN_\f^*\rN_\f={\rm span}\{y^*x\;:\; x,\, y\in \rN_\f\}\,.$$
It is well-known that  for any $x\in\rM_\f$ the operator
$D^{1/2}xD^{1/2}$ is closable and its closure belongs to $L_1(\M,
\f)$ (cf. e.g. \cite{terp-int} and \cite{gold-linKMSS}). Denoting
the closure again by $D^{1/2}xD^{1/2}$, we then have
  \be
  \f(x)=\tr(D^{1/2}xD^{1/2}),\quad x\in\rM_\f\,.
  \ee
In particular, if $\f$ is bounded, then $D\in L_1(\M,\f)$ and
 \beq\label{state vs trace}
 \f(x)=\tr(D^{1/2}xD^{1/2})=\tr(Dx),\quad x \in\M.
 \eeq

\begin{rk}\label{subLp}
 Let $\N$ be a w*-closed involutive subalgebra of $\M$ such that
$\f\big|_{\N}$  is semifinite. Assume that $\N$ is invariant under
$\s^\f$, i.e., $\s_t^\f(\N)\subset \N$ for all $t\in\real$. Then
it is easy to check that $L_p(\N,\f\big|_{\N})$ coincides
isometrically with a subspace of $L_p(\M,\f)$ for every $0<p<\8$.
 \end{rk}

\begin{rk}\label{uniqueness Lp}
 It is proved in \cite{haag-Lp} and \cite{terp} that
$L_p(\M, \f)$ is independent of $\f$ up to isometric isomorphism
preserving the order and modular structure of $L_p(\M, \f)$. This
independence allows us to denote $L_p(\M,\f)$ simply by $L_p(\M)$.
On the other hand, if $\f$ is tracial, i.e., $\f(x^*x)=\f(xx^*)$
for all $x\in\M$, then the Haagerup space $L_p(\M, \f)$
isometrically coincides with the tracial noncommutative
$L_p$-space associated with $\f$ constructed by Dixmier
\cite{dixLp} and Segal \cite{segalLp}.  We refer to
\cite{px-survey} for more information and historical references on
noncommutative $L_p$-spaces.
 \end{rk}


\section{The reduction theorem}
\label{Reduction}


This section is the core of the paper. Its result is a reduction
theorem that approximates a type III von Neumann algebra by finite
ones. This theorem is due to the first-named author and has never
appeared in published form. Our presentation will follow the
first-named author's manuscript, which has been circulated in a
limited circle of people since \cite{haag-red}. We will
concentrate our attention mainly on the $\sigma$-finite case.
Throughout this section, $G$ will denote the discrete subgroup
$\bigcup_{n \geq 1} 2^{-n} \ent$ of $\real$. Let $\M$ be a
$\s$-finite von Neumann algebra equipped with a normal faithful
state $\f$. Consider the crossed product $\R = \M
\rtimes_{\s^{\f}} G$. Here, the modular automorphism group
$\s^{\f}$ is also viewed as an automorphic representation of $G$
on $\M$. Let $\wh\f$ denote the dual weight of $\f$. Recall that
$\wh\f$ is a normal faithful state on $\R$.

 \begin{thm}\label{red s-finite}
 With the notation above, there exists an increasing sequence
$(\R_{n})_{n \geq 1}$ of von Neumann subalgebras of $\R$
satisfying the following properties:
 \begin{enumerate}[\rm(i)]
 \item each $\R_{n}$ is finite;
 \item $\bigcup_{n\geq 1}\R_n$ is w*-dense in $\R$;
 \item for every $n\in\nat$ there exists a normal faithful
conditional expectation
 $\Phi_{n}$ from $\R$ onto $\R_{n}$ such that
  $$\wh\f\circ\Phi_{n}=\wh\f \quad
  \mbox{and } \quad
  \s_{t}^{\wh\f}\circ\Phi_{n}=\Phi_{n}\circ
  \s^{\wh\f}_{t}, \quad t\in\real.$$
  \end{enumerate}
  \end{thm}

The remainder of this section is mainly devoted to the proof of
this theorem. $\R_{n}$ will be constructed as the centralizer of a
normal faithful positive functional $\f_{n}$ such that the modular
automorphism group $\s_{t}^{\f_{n}}$ is periodic with period
$2^{-n}.$ We  keep all the notation introduced in
subsection~\ref{Crossed products}. $\T$ denotes the unit circle of
the complex plane equipped with normalized Haar measure $dm$.

\begin{lem}\label{distribution of lambda}
 For any $f \in L_\8(\T)$,
 \beq\label{distribution of lambda formula}
 \wh\f(f(\l(t))) = \int_{\T}\ f(z)dm(z),\ \
 t\in G \setminus\{ 0 \}.
 \eeq
In other words, the distribution measure of $\l(t)$ with respect
to $\wh\f$ is equal to $dm$.
 \end{lem}

\pf Let $t \in G\setminus\{ 0 \}$. Then by \eqref{dual state}, for
any $n \in \ent$,
 $$\wh\f\big(\l (t)^{n}\big) =
 \wh\f(\l (nt)) = \left\{\begin{array}{ll}
  1 & \textrm{ if } n=0, \\
  0  &\textrm{ if } n\neq 0.
  \end{array}\right.$$
Thus \eqref{distribution of lambda formula} holds whenever $f$ is
a monomial $z^{n}$, $n\in\ent$, so also whenever $f$ is a
trigonometric polynomial. Then the normality of $\wh\f$ yields
\eqref{distribution of lambda formula} for all $f \in L_\8(\T)$.
\cqd

\smallskip

Recall that if $\psi$ is an {\it n.s.f.} weight on a von Neumann
algebra $\N$, the centralizer $\N_{\psi}$ of $\psi$ is the fixed
point algebra of $\s^{\psi}_{t}$:
 $$\N_{\psi} = \{x\in\N\; :\; \s^{\psi}_{t}(x) = x,\
 \forall\; t \in\real\}.$$
If additionally $\psi$ is bounded,
 $$\N_{\psi} = \{ x\in \N\; :\; \psi x = x \psi\},$$
where for $a, b\in\N$, $a \psi b$ denotes the functional on $\N$
given by
 $$a\psi b(y) = \psi (bya),\ y \in \N.$$
$\Z(\N)$ denotes the center of $\N.$

\begin{lem}\label{bn}
 \begin{enumerate}[\rm(i)]
 \item $\l(t) \in \Z(\R_{\wh\f})$ for any $t\in G$.
 \item For every $n \in\nat$ there exists a unique $b_{n} \in
\Z(\R_{\wh\f})$ such that $0 \leq b_{n} \leq 2 \pi$ and $e^{ib_n}
= \l(2^{-n})$.
 \end{enumerate}
 \end{lem}

\pf By \eqref{com pi-l} and \eqref{modular gp of dual weight},
 \beq\label{modular gp of dual weight bis}
 \s_{s}^{\wh\f} (x) = \l(s)x
 \l(s)^{\ast},\ x \in \R, \ s \in G.
 \eeq
This clearly implies (i). To prove (ii) we use the principal
branch ${\rm Log}\,z$ of the logarithmic function which satisfies
 $0\leq{\rm Im}({\rm Log}\,z)<2\pi,\ z \in\com\setminus\{0\}$.
Let
 $$b_{n} = -i {\rm Log}\big(\l (2^{-n})\big).$$
Then $0 \leq b_{n} \leq 2 \pi,\ e^{i b_n} = \l (2^{-n})$; by (i)
and functional calculus, $b_{n} \in \Z(\R_{\wh \f}).$ The
uniqueness of $b_{n}$ follows from the fact that $\l(2^{-n})$ has
no point spectrum by virtue of Lemma \ref{distribution of lambda}
and the faithfulness of $\wh\f.$ \cqd

\smallskip

Now let $a_{n} = 2^{n} b_{n},$ and define a sequence $(\f_{n})_{n
\geq 1}$ of normal faithful positive functionals on $\R$  by
 \beq\label{def of fn}
 \f_{n}(x)=\wh\f(e^{-a_{n}} x),\quad
 x \in \R,\ n \ge 1.
 \eeq

\begin{lem}\label{construction of Rn}
 \begin{enumerate}[\rm(i)]
 \item $\s_{t}^{\f_{n}}$ is $2^{-n}$-periodic for all $n\ge1$.
 \item Let $\R_{n} = \R_{\f_{n}}, \ n \geq 1.$
There exists a unique normal faithful conditional expectation
$\Phi_{n}$ from $\R$ onto $\R_{n}$ such that
 $$\wh\f\circ\Phi_{n} = \wh\f
 \quad\mbox{and }\quad
 \s_{t}^{\wh \f}\circ\Phi_{n} = \Phi_{n}\circ
 \s_{t}^{\wh \f}, \quad t \in\real,\ n\ge1.$$
 \item $\R_{n} \subset \R_{n+1}$.
 \end{enumerate}
 \end{lem}

\pf (i) By Lemma \ref{bn} (i), $a_{n} \in \Z(\R_{\wh\f})$; in
particular, $a_{n} \in \R_{\wh\f}.$ Thus
 \beq\label{formula od sfn}
 \s_{t}^{\f_{n}}(x)=e^{-ita_{n}}\s_{t}^{\wh\f} (x)e^{ita_{n}},
 \quad x \in\R,\ t \in\real.
 \eeq
Therefore, by \eqref{modular gp of dual weight bis} and Lemma
\ref{bn} (ii), for all $x\in\R,$
 $$\s_{2^{-n}}^{\f_{n}}(x)
 = e^{-i\,b_{n}}\s^{\wh \f}_{2^{-n}}(x)e^{ib_{n}}
 = \l(2^{-n})^{\ast}\l(2^{-n})x
 \l(2^{-n})^{\ast}\l(2^{-n})
 = x,$$
whence (i).

(ii)  Define $\Phi_{n}$ by
 $$\Phi_{n}(x) = 2^{n}\int^{2^{-n}}_{0} \s_{t}^{\f_{n}}
 (x) dt,\quad x \in \R.$$
By  the $2^{-n}$-periodicity of $\s^{\f_{n}}$, we find
 $$\Phi_{n}(x) = \int^{1}_{0}\s_{t}^{\f_{n}}(x)dt,
 \quad x \in \R.$$
It is then a routine matter to check that $\Phi_{n}$ is a normal
faithful conditional expectation from $\R$ onto $\R_{n}.$ Since
$a_{n} \in \R_{\wh \f},$ by \eqref{formula od sfn} and the
$\s^{\wh\f}$-invariance of $\wh\f,$ we get
 $$\wh\f\big(\s_{t}^{\f_{n}} (x)\big)
 =\wh\f \big(\s_{t}^{\wh \f}(x) \big)
 =\wh\f(x),\quad x \in \R,\ t \in\real.$$
Thus $\wh\f$ is also $\s^{\f_{n}}$-invariant for all $n \geq 1.$
It then follows that
 $$\wh\f(\Phi_{n}(x))=\int^{1}_{0}\wh\f\big(\s_{t}^{\f_{n}}(x)\big)dt
 =\wh\f (x),\quad x \in \R,$$
so $\wh\f\circ \Phi_{n}=\wh \f.$ This latter equality implies the
uniqueness of $\Phi_{n}$ as well as the commutation relation
between $\Phi_{n}$ and $\s_{t}^{\wh \f}$ (cf., e.g.,
\cite{tak-cond}). Alternately, this commutation relation
immediately follows from \eqref{formula od sfn} and the definition
of $\Phi_{n}.$

(iii) For every natural number $n$, $a_n$ and $a_{n+1}$ commute,
because they are both contained in $\Z(\R_{\wh\f})$. Hence by
\eqref{formula od sfn}, $\f_{n+1}$ is $\s^{\f_n}$-invariant. Let
$h_{n} = \displaystyle \frac{D\f_{n+1}}{D\f_{n}},$ the
Radon-Nikod\'ym derivative of $\f_{n+1}$ with respect to $\f_{n}$
(cf. \cite{ped-tak}). Then $\f_{n+1}(x) = \f_{n}(h_{n} x)$ for all
$x\in\R$, and by the definition of $\f_{n}$,
 $$h_{n} = \frac{D\f_{n+1}}{D\wh\f}\cdot \frac{D\wh\f}{D\f_{n}}
 =e^{-(a_{n+1} - a_{n})}.$$
Proving $\R_{n}\subset\R_{n+1}$  is equivalent to showing $h_{n}
\in \Z(\R_{n})$ for every $n\ge1$. By \eqref{formula od sfn},
$\R_{\wh \f} \subset \R_{\f_{n}}$ since $a_{n}\in \Z(\R_{\wh\f})$.
In particular, $h_{n} \in \R_{\f_{n}}=\R_n$. Now
 $$a_{n} = -i2^{n}{\rm Log}\,\l (2^{-n}) = -i2^{n}{\rm
 Log}\big(\l(2^{-n-1})^{2}\big).$$
Thus
 $$a_{n+1} - a_{n} = -i2^{n} \big[ 2{\rm Log}\,
 \l(2^{-n-1}) -{\rm Log}\,\big(\l
 (2^{-n-1})^{2}\big)\big].$$
However, for any $z \in\T$,
 $$2\,{\rm Log}\, z - {\rm Log}(z^{2}) = \left\{\begin{array}{ll}
 0 & \textrm{if}\  0 \le{\rm Arg}\,z <\pi,\\
 2 \pi i &\textrm{if}\ \pi \le {\rm Arg}\,z <2 \pi.
 \end{array}\right.$$
Hence
 $$a_{n+1} - a_{n} = 2^{n+1}\pi  e_{n},$$
where $e_{n}$ is the spectral projection of $\l (2^{-n-1})$
corresponding to ${\rm Im}\, z < 0.$ Therefore, for all $x\in\R$
and $t\in\real$,
 $$\s_{t}^{\f_{n+1}}(x) = h_{n}^{it}\s_{t}^{\f_{n}}(x)h_{n}^{-it}
 = e^{-i 2^{n+1} \pi te_{n}}
 \s_{t}^{\f_{n}} (x) e^{i 2^{n+1} \pi te_{n}}\,.$$
Consequently, if $x \in \R_n$, then by the $2^{-n-1}$-periodicity
of $\s_{t}^{\f_{n+1}}$, we deduce
 $$x = e^{-i \pi e_{n}} x e^{i \pi e_{n}} = (1-2 e_{n}) x(1-2e_{n}).$$
Therefore, $1-2e_{n} \in \Z(\R_n)$,  so $e_{n} \in \Z(\R_n)$. Thus
$h_{n} \in \Z(\R_n)$, which in turn yields the desired inclusion
$\R_{n} \subset \R_{n+1}$. \cqd

\smallskip

By Lemma \ref{construction of Rn}, in order to complete the proof
of Theorem \ref{red s-finite}, it remains to show the w*-density
of the union of the $\R_{n}$ in $\R$. This is the most difficult
part of the proof, and will be done in the following lemmas. The
first one is well known. Recall that an element $x\in\N$ is called
analytic (with respect to $\s^\psi$) if the map $t\mapsto
\s^\psi_t(x)$ extends to an entire function, i.e., if there exists
a (necessarily unique) analytic function $F_x: \com\to \M$ such
that $F_x(t) = \s^\psi_t(x)$ for all $t\in\real$. In this case we
put $\s^\psi_s(x) = F_x(s)$ for all $s\in\com$. Let $\N_a$ denote
the family of all analytic elements of $\N$.

\begin{lem}\label{L2 bounded elements}
 Let $\N$ be a von Neumann algebra and $\psi$ a normal
faithful state on $\N$.
 \begin{enumerate}[\rm(i)]
 \item If $x\in\N_a$, then there exists a constant $c\geq 0$ such that
$x \psi x^{\ast} \leq c \psi.$
 \item If $x \in \N_{\psi},$ then $x \psi x^{\ast} \leq \| x
\|^{2} \psi.$
 \end{enumerate}
 \end{lem}

\pf Without loss of generality, we assume that $\N$ acts
standardly on a Hilbert space $H.$ Then there exists a cyclic and
separating vector $\xi_{0} \in H$ such that
 $\psi(x) = \la x\xi_{0},\ \xi_{0}\ra.$
Let $\Delta$ and $J$ be the corresponding modular operator and
isometric conjugation (see, for instance,  \cite[section
9.2]{kar-II}). Then
 $$\s_{t}^{\psi} (x) = \Delta^{it}x\Delta^{-it},\ x \in
 \N,\ t \in {\real}.$$
Let $x\in\N_a$. Then
 \be
 x\xi_{0}
 &=& J\Delta^{1/2}x^{\ast}\xi_{0}
 = J\Delta^{1/2}x^{\ast}\Delta^{-1/2}\xi_{0}\\
 &=& J\s^{\psi}_{{-i /2}} (x^{\ast})\xi_{0}
 = \big(J\s^{\psi}_{{-i /2}}(x^{\ast})J \big)(\xi_{0}).
 \ee
Now $x' = J\s^{\psi}_{{-i /2}} (x^{\ast}) J\in \N'$. Hence, for
any $y \in \N$ and $y\geq 0$,
 $$\psi (x^{\ast} y x ) = \langle yx'\xi_{0},\ x' \xi_{0} \rangle
 \leq \| x' \|^{2} \langle y\xi_{0},\ \xi_{0}\rangle,$$
whence $x \psi x^{\ast} \leq \| x' \|^{2} \psi.$ If additionally
$x \in \N_{\psi},$ then $\s_{{-i/2}} (x^{\ast}) = x^{\ast}$, so
$x' = J x^{\ast} J.$ Thus (ii) follows. \cqd

\smallskip

In the following lemma, $[x,\;y]$ denotes the commutator of two
operators $x$ and $y$, i.e., $[x,\; y]=xy-yx$. If $\psi$  is a
positive functional on an algebra $\N$,  $\|x\|_\psi^2=\psi(x^*x)$
for $x\in\N$.

\begin{lem}\label{L2 estimate}
 Keep the notation in Lemma \ref{bn}. Then for $x\in\R$,
 \begin{enumerate}[\rm(i)]
  \item $\displaystyle\lim_{n \to\infty}
  \| [b_{n},\; x]\|_{\wh\f} = 0;$
  \item $\displaystyle\lim_{n \to\infty}\,\sup_{t\in
 [-1,\; 1]}\ \bigl\|\big[e^{i t b_{n}},\; x \big]\big\|_{\wh\f}=0;$
 \item $\displaystyle\lim_{n\to\infty}\,\sup_{t\in\real} \,
 \big\|\s_{t}^{\wh\f}(x) - e^{i a_{n} t} xe^{-i a_{n}t}\big\|_{\wh\f}=0.$
 \end{enumerate}
 \end{lem}

\pf (i) Let $x \in \R$ and $k\in\ent$. Since $\l(t)\in\R_{\wh\f}$
because of Lemma \ref{bn}, we have
 \be
 \big\|[\l (2^{-n})^{k},\; x]\big\|_{\wh\f}
 &=& \big\|\l(k2^{-n})x - x \l(k2^{-n})\big\|_{\wh\f}\\
 &=& \big\|\l (k2^{-n})x \l(k2^{-n})^{\ast} - x\big\|_{\wh\f}
 = \big\|\s_{k2^{-n}}^{\wh\f}(x)- x\big\|_{\wh\f}\,.
 \ee
It follows that
 \beq\label{L2 estimate 1}
 \lim_{n\to\infty}\big\|[ P(\l(2^{-n})),\;
 x]\big\|_{\wh\f} = 0
 \eeq
for any monomial, so for any trigonometric polynomial $P$.

Now assume that $x$ is analytic with respect to $\s_{t}^{\wh\f}.$
Since ${\rm Log}\in L_2(\T),$  given $\e> 0$ there exists a
trigonometric polynomial $P$ such that
 $$\| P + i\ {\rm Log}\|_{L_2(\T)} < \e.$$
Recalling that $b_{n} = -i{\rm Log}\big(\l (2^{-n})\big),$ we
deduce from \eqref{distribution of lambda formula} that
 $$\|b_{n} - P(\l (2^{-n}))\|_{\wh\f}
 =\big(\wh\f(| b_{n} -P(\l(2^{-n}))|^{2})\big)^{1/2}
 =\|-i{\rm Log}\, - P \|_{L_2(\T)}<\e.$$
On the other hand, since $x$ is analytic, by Lemma \ref{L2 bounded
elements} there exists a constant $c$ such that for any $y \in
\R$,
 $$\| [y,\; x ]\|_{\wh\f}\le\|yx\|_{\wh\f}+\| xy \|_{\wh\f}
 \le\sqrt{c}\,\| y \|_{\wh\f} + \|x\| \, \|y\|_{\wh\f}
 = (\sqrt{c} + \|x\|)\|y\|_{\wh\f}\,.$$
Combining the preceding inequalities, we obtain
 \be
 \| [b_{n},\; x ]\|_{\wh\f}
 &\le& \|[ P(\l(2^{-n})),\; x]\|_{\wh\f}
 + \| [b_{n} -P(\l(2^{-n})),\; x]\|_{\wh\f}\\
 &\le& \| [P(\l(2^{-n})),\; x]\|_{\wh\f}
 + (\sqrt{c} + \|x\|)\|b_{n}-P(\l (2^{-n}))\|_{\wh\f}\\
 &\le& \|[P(\l (2^{-n})),\; x]\|_{\wh\f}
 + \e(\sqrt{c} + \| x\|).
 \ee
Therefore, by \eqref{L2 estimate 1},
 $$\limsup_{n\to\infty}\|[b_{n},\; x]\|_{\wh\f}\le\e(\sqrt{c}+\|x\|),$$
whence
 $$\lim_{n \to\infty}\| [b_{n},\; x]\|_{\wh\f} = 0, \quad
 \forall\; x\in\R_a\,,$$
where $\R_a$ denotes the set of analytic elements in $\R$ with
respect to $\s^{\wh\f}$.

Next, fix any $x\in\R$. Since  $\R_a$ is a $\s$-strongly dense
involutive subalgebra of $\R$, for every $\e>0$, we can choose an
$x_{0}\in\R_a$ such that $\|x-x_{0}\|_{\wh\f}<\e$. Then by the
fact that $b_{n}\in\R_{\wh\f}$ and $\| b_{n} \| \leq 2 \pi$ from
Lemma \ref{bn}, we deduce, as before,
 $$\lim_{n \to \infty}\|[b_{n},\; x]\|_{\wh\f} =0.$$

(ii) For any $k\in\nat,$ we have
 $$[b_{n}^{k},\; x]=b_{n}[b_{n}^{k-1},\; x]+[ b_{n},\; x]b_{n}^{k-1}.$$
Since $b_{n} \in \R_{\wh \f}\,$, an induction argument yields
 $$\| [ b_{n}^{k},\; x]\|_{\wh \f}
 \leq k \| b_{n} \|^{k-1}\| [b_{n},\ x]\|_{\wh\f}
 \le k(2\pi)^{k-1}  \| [ b_{n},\;x]\|_{\wh\f}.$$
Hence for any $z \in \com$,
 \be
 \| [e^{z b_{n}},\; x]\|_{\wh\f}
 &\le& \sum_{k\ge1}\frac{|z|^{k}}{k!}\,\|[ b^{k}_{n},\;
 x]\|_{\wh\f}\\
 &\le& \sum_{k\ge1}\frac{|z|^{k}}{(k-1)!}\,
  (2\pi)^{k-1} \|[b_{n},\; x]\|_{\wh\f}\\
  &=&|z|\,  e^{2\pi|z|}\|[ b_{n},\; x]\|_{\wh\f}.
  \ee
Therefore
 $$\sup_{t\in [-1,\; 1]}\|[ e^{it b_{n}},\; x]\|_{\wh\f}
 \le e^{2 \pi} \|[ b_{n}, \; x]\|_{\wh\f},$$
which, together with (i), implies (ii).

(iii) Fix $x \in \R$ and $\e> 0.$ By (ii) there exists $n_{0} \in
\nat$ such that
 \beq\label{L2 estimate2}
 \|[ e^{i s b_{n}},\; x]\|_{\wh\f}\le\e,
 \quad \forall\;s \in [-1,\; 1],\ \forall\;n\ge n_{0}.
 \eeq
Moreover, $n_0$ can be chosen such that further
 \beq\label{L2 estimate3}
 \|\s_{s}^{\wh\f}(x) - x \|_{\wh\f}\le \e,\quad |s| \le 2^{-n_{0}}.
 \eeq
Let $t\in\real$ and $n\in\nat$ with $n\ge n_{0}.$ Write $t = t_{1}
+ t_{2},$ where $t_{1} = k\, 2^{-n}$ for some $k\in\ent$ and $0
\le t_{2} \le 2^{-n}.$ Then for any $y \in \R$,
 \be
 \s_{t_{1}}^{\wh\f} (y)
 &=& \l(k 2^{-n})y\l(k2^{-n})^{\ast}
 =e^{i k b_{n}}y e^{-ik b_{n}}\\
 &=& e^{i k 2^{-n} a_{n}}ye^{-ik2^{-n}a_{n}}
 = e^{i t_{1}a_{n}}ye^{-i t_{1}a_{n}}.
 \ee
Since $\|\cdot\|_{\wh\f}$ is invariant under $\s_{t_{1}}^{\wh\f}$
and $a_{n} \in \Z(\R_{\wh\f}),$ we deduce
 \be
 \| \s_{t}^{\wh \f}(x) - e^{i a_{n}t} x
 e^{-ia_{n}t}\|_{\wh\f}
 &=& \| \s_{t_{2}}^{\wh\f} (x) - e^{i a_{n} t_{2}}
 x e^{-i a_{n}t_{2}}\|_{\wh\f}\\
 &\le& \|\s_{t_{2}}^{\wh \f} (x) - x \|_{\wh\f}
 + \| x - e^{i a_{n}t_{2}}x e^{-ia_{n}t_{2}}
 \|_{\wh \f}\\
 &=& \|\s_{t_{2}}^{\wh \f}(x) - x \|_{\wh\f}
 + \|[e^{-ia_{n}t_{2}},\; x]\|_{\wh\f}\,.
 \ee
Now $a_{n}t_{2}=(2^{n} t_{2}) b_{n}$ and $2^{n}t_{2}\le 1$. Hence
from \eqref{L2 estimate2} and \eqref{L2 estimate3} it follows that
 $$\| \s_{t}^{\wh\f}(x)
 -e^{ia_{n}t}xe^{-ia_{n} t}\|_{\wh\f}\le 2\e.$$
This yields (iii). \cqd

\smallskip

Finally, we are ready to show the w*-density of $\bigcup_{n \geq
1}\ \R_{n}$ in $\R$.

\begin{lem}
 For any $x \in \R,$ $\Phi_{n}(x)$ converges to $x$
$\s$-strongly as $n\to\infty.$ Consequently, $\bigcup_{n\ge1}\
\R_{n}$ is $\s$-strongly dense  in $\R.$
 \end{lem}

\pf By the definition of $\Phi_{n},$ it suffices to show
 $$\lim_{n\to\infty}\sup_{t\in\real}
 \big\|\s_{t}^{\f_{n}} (x) - x \big\|_{\wh \f}
 = 0,\quad x\in \R.$$
By \eqref{formula od sfn} and the fact that $a_{n} \in
\R_{\wh\f}$, we have
 $$
 \big\|\s_{t}^{\f_{n}}(x)-x\big\|_{\wh\f}
 =\big\|e^{-ita_{n}}\s_{t}^{\wh\f}(x)e^{ita _{n}}-x\big\|_{\wh\f}
 =\big\|\s_{t}^{\wh \f}(x) - e^{ita_{n}}xe^{-ita_{n}}\big\|_{\wh\f}\,.$$
Therefore, the desired limit follows from Lemma \ref{L2 estimate}
(iii). Thus the proof of Theorem 2.1 is complete.
 \cqd

\smallskip

Theorem \ref{red s-finite} can be extended to general von Neumann
algebras (not necessarily $\s$-finite) as follows.

\begin{rk}\label{red general}
  Let $\M$ be a von Neumann algebra.
Then there exist a von Neumann algebra $\R$ and an increasing
family $\{\R_{i}\}_{i\in I}$ of w*-closed involutive subalgebras
of $\R$ satisfying the following properties:
 \begin{enumerate}[\rm(i)]
 \item $\M$ is a von Neumann subalgebra of $\R$ and there exists
a normal faithful conditional expectation $\Phi$ from $\R$ onto
$\M$;
 \item each $\R_{i}$ admits a normal faithful tracial state;
 \item the union of all $\R_{i}$ is w*-dense in $\R$;
 \item for every $i\in I$ there exists a normal conditional
expectation $\Phi_{i}$
 from $\R$ onto $\R_{i}$ such that
 $$\Phi_{i}\circ \Phi_{j}=\Phi_{j}\circ\Phi_{i}
 =\Phi_{i}\quad\mbox{whenever}\quad i\le j;$$
 \item there exists an {\it n.s.f.} weight  $\f$ on
$\M$ such that
 $$\wh\f\circ\Phi_{i}=p_i\wh\f p_{i} \quad\mbox{and}\quad
 \s_{t}^{\wh \f}\circ\Phi_{i}
 =\Phi_{i}\circ\s_{t}^{\wh \f},\quad t\in\real,\ i\in I,$$
where $\wh\f=\f\circ\Phi$ and $p_i$ is the identity of $\R_i$.
 \end{enumerate}
 \end{rk}

It is not hard to deduce this statement from Theorem~\ref{red
s-finite}. Indeed, this is quite easy if $\M=\N\bar\ot B(K)$ for
some $\s$-finite von Neumann algebra $\N$ and some Hilbert space
$K$. The general case can be reduced to the previous one by using
the classical fact that any von Neumann algebra is a direct sum of
algebras of the form $\N\bar\ot B(K)$ with $\N$ $\s$-finite.
Namely,  any $\M$ admits the following direct sum decomposition:
 $$\M=\bigoplus_{j\in J}\N_j\bar\ot B(K_j),$$
where each  $\N_j$ is a $\s$-finite von Neumann algebra. See, for
instance, the proof of Theorem~2.8.1 in Sakai's book \cite{sakai}.


\section{Applications to noncommutative $L_p$-spaces}
  \label{Applications to noncommutative Lp-spaces}


In this section we use Theorem~\ref{red s-finite} to prove an
approximation theorem of general noncommutative $L_p$-spaces by
those associated with finite von Neumann algebras. This is due to
the first-named author and is indeed the original intention of
\cite{haag-red}.

\begin{thm}\label{approximation s-finite}
 Let $\M$ be a $\s$-finite von Neumann algebra and
$0 < p < \infty.$ Let $L_p(\M)$ be the Haagerup noncommutative
$L_p$-space associated with $\M$. Then there exist a Banach space
$X_p$ $($a quasi-Banach space if $p<1),$ a sequence
$(\R_n)_{n\ge1}$ of finite von Neumann algebras, each equipped
with a normal faithful finite trace $\t_n$, and for each $n\ge1$
an isometric embedding $J_n\, :\, L_p(\R_n,\t_n) \to X_p$ such
that
 \begin{enumerate}[\rm(i)]
  \item the sequence $\big(J_n\big(L_p(\R_n, \t_n)
 \big)\big)_{n\ge1}$ is increasing;
 \item  $\bigcup_{n\ge1}J_n\big(L_p(\R_n,\t_n)\big)$
is dense in $X_p$;
 \item  $L_p(\M)$ is isometric to a subspace $Y_p$ of $X_p$;
 \item  $Y_p$ and all $J_n\big(L_p(\R_n, \t_n)\big)$ are
$1$-complemented in $X_p$ for $1\le p <\8$.
 \end{enumerate}
Here $L_p(\R_n, \t_n)$ is the tracial noncommutative $L_p$-space
associated with $(\R_n,\t_n)$.
 \end{thm}

\pf Fix a normal faithful state $\f$ on $\M$. We will use
Theorem~\ref{red s-finite} and keep all the notation there. The
space $X_p$ required in the statement above will be $L_p(\R,
\wh\f)$. By Remark~\ref{subLp}, $L_p (\M, \f)$ and all $L_p(\R_n,
\wh\f\big|_{\R_n})$ are naturally isometrically identified as
subspaces of $L_p(\R, \wh \f)$ for $0<p<\8$. Moreover, the
sequence $\big(L_p(\R_n, \wh\f\big|_{\R_n})\big)_{n\ge1}$ is
increasing. On the other hand, by \cite[Lemma 2.2]{ju-doob},
$\bigcup_{n\ge1}\, L_p(\R_n, \wh\f_n)$ is dense in $L_p(\R,
\wh\f)$ for $0 < p < \infty$. Finally, since each $\R_n$ is a
finite von Neumann algebra with a finite normal faithful trace
$\t_n$, $L_p(\R_n, \wh\f_n)$ is isometric to the usual
noncommutative $L_p$-space on $\R_n$ defined by $\t_n$ (see
Remark~\ref{uniqueness Lp}). Hence, the space $L_p(\R, \wh\f)$ and
the sequence $\big(L_p(\R_n, \wh\f_n)\big)_{n\ge1}$ satisfy
properties (i) - (iii). The complementation property for $p\ge1$
in (iv) follows from \cite[Lemma~2.2]{jx-burk} thanks to the
conditional expectations $\Phi$ and $\Phi_n$. \cqd

\smallskip

The following two remarks show that Theorem \ref{approximation
s-finite} is general enough for most applications.

\begin{rk}\label{s-finite-general}
 Let $\M$ be a general von Neumann algebra. Then for any $0<p<\8$,
 $$L_p(\M)=\bigcup_e eL_p(\M)e=\bigcup_e L_p(e\M e),$$
where the union runs over the directed set of all $\s$-finite
projections of $\M$. Indeed, for any $x\in L_p(\M)$, the left and
right support projections of $x$ are $\s$-finite, so is their
union $e$. Thus $x\in eL_p(\M)e$.
 \end{rk}

\begin{rk}
 Theorem \ref{approximation s-finite}
(combined with Remark~\ref{s-finite-general}) reduces many
geometrical properties of general noncommutative $L_p$-spaces to
the corresponding ones in the tracial case. This is indeed true
for all those properties which are of a finite-dimensional nature.
These include, for instance, Clarkson's inequalities, uniform
convexity, uniform smoothness, type, cotype, and UMD property. We
refer to \cite{px-survey} for more information.
 \end{rk}

\begin{rk}\label{approximation}
 Using Remark~\ref{red general}, we can extend
Theorem \ref{approximation s-finite} to the general case as
follows.  Let $\M$ be a general von Neumann algebra and $0 < p <
\infty.$ Let $L_p(\M)$ be the Haagerup noncommutative $L_p$-space
associated with $\M$. Then there exist a Banach space $X_p$ $($a
quasi-Banach space if $p<1),$ a family $\{\R_{i}\}_{i \in I}$ of
finite von Neumann algebras, each equipped with a normal faithful
finite trace $\t_{i}$, and for each $i \in I$ an isometric
embedding $J_{i}\, :\, L_p(\R_{i},\t_{i}) \to X_p$ such that
 \begin{enumerate}[\rm(i)]
  \item $J_{i}\big(L_p(\R_{i}, \t_{i})\big)\subset
  J_{j}\big(L_p(\R_{j}, \tau_{j})\big)$ for all $i, j \in I$
such that $i \leq j$;
 \item  $\bigcup_{i\in I}J_{i}\big(L_p(\R_{i},\t_{i})\big)$
is dense in $X_p$;
 \item  $L_p(\M)$ is isometric to a subspace $Y_p$ of $X_p$;
 \item  $Y_p$ and $J_{i}\big(L_p(\R_{i}, \t_{i})\big)$, $i\in I$ are
$1$-complemented in $X_p$ for $1\le p <\8$.
 \end{enumerate}
 \end{rk}


\section{Extensions of maps to crossed products}
 \label{Extensions of maps to crossed products}


In noncommutative analysis we often need to extend a map between
two von Neumann algebras to their corresponding noncommutative
$L_p$-spaces. On the other hand, when applying Theorem \ref{red
s-finite} to concrete problems we also need to extend maps between
von Neumann algebras to their crossed products. This section is
devoted to the second type of extensions, while the next section
is devoted to the first type. In  what follows, all maps
considered will be assumed linear.

\begin{thm}\label{extension crossed}
 Let $\M$ and $\N$ be two von Neumann algebras acting on the
same Hilbert space $H,$ $G$ a locally compact abelian group, $\a$
and $\b$ two automorphic representations of $G$ on $\M$ and $\N$,
respectively. Assume that $T\, :\, \M \to \N$ is a completely
bounded normal map such that
 \beq\label{Tab}
 T\circ\a_{g}=\b_{g}\circ T,\quad g \in G.
 \eeq
Then $T$ admits a unique completely bounded normal extension $\wh
T$ from $\M\rtimes_{\a} G$ into $\N\rtimes_{\b} G$ such that
$\|\wh T\|_{cb}=\|T\|_{cb}$ and
 \beq\label{def of T hat}
 \wh T\big(\l(g) \pi_{\a}(x)\big) = \l(g)
 \pi_{\b}(T x),\quad x \in \M,\  g \in G.
 \eeq
Moreover, $\wh T$ satisfies the following properties:
 \begin{enumerate}[\rm(i)]
 \item Let $\A$ be the von Neumann subalgebra on $L_2(G, H)$
generated by all $\l(g),\,  g \in G.$ Then
 \beq\label{def of T hat1}
 \wh T\big(a\pi_{\a}(x)b\big) = a\pi_{\b}(T x)b,\quad
 x \in \M,\ a,\; b \in \A.
 \eeq
 \item $\wh T\circ\wh\a_{\g}=\wh\b_{\g}\circ\wh T$, $\g\in\wh G.$
 \item If $T$ is a homomorphism, $\ast$-homomorphism or
completely positive map, so is $\wh T.$
 \item If $\N$ is a subalgebra of $\M$, $\b=\a\big|_\N$ and
$T$ is a $($faithful$)$ normal conditional expectation from $\M$
onto $\N$, then $\wh T$ is a $($faithful$)$ normal conditional
expectation from $\M\rtimes_{\a} G$ onto $\N\rtimes_{\b} G$.
 \item Let $\f$ $($resp. $\psi)$ be an {\it n.s.f.} weight
on $\M$ $($resp. $\N)$ such that
 \beq\label{Tsf}
 T\circ \s^{\f}_{t} = \s_{t}^{\psi} \circ T,\quad  t \in
 \real.
 \eeq
Then
 \beq\label{That sfhat}
 \wh T \circ \s^{\wh \f}_{t} =
 \s_{t}^{\wh \psi} \circ \wh T,\quad t \in\real.
 \eeq
 \item Assume in addition that $T\ge0$. Let $\f$ $($resp. $\psi)$
be an {\it n.s.f.} weight on $\M$ $($resp. $\N)$ such that
$\psi\circ T\le\f$. Then $\wh\psi\circ\wh T\le\wh\f$.

 \end{enumerate}
 \end{thm}

\pf Since $T$ is completely bounded and normal, so is
 $$T \ot {\rm id}_{B(L_2(G))}\,:\, \M\bar\ot B(L_2(G))
 \to \N\bar\ot B(L_2(G)).$$
Moreover,
 $$\|T \ot {\rm id}_{B(L_2(G))}\|_{cb}
 = \|T\|_{cb}\,.$$
We claim that $T \otimes {\rm id}_{B(L_2(G))}$ maps $\M
\rtimes_{\a} G$ into $\N \rtimes_{\b} G.$ Fix an orthonormal basis
$(f_i)_{i\in I}$ in $L_2(G)$. Then by the definition of $\pi_{\a}
(x),$ one sees that the matrix of $\pi_{\a} (x)$ in this basis has
its coefficient at the position $(i, j)$ given by
 $$\big(\pi_{\a}(x)\big)_{ij}
  =\int_{G}\a^{-1}_{h}(x)\bar f_{i}(h)f_{j}(h)dh.$$
Thus by the normality of $T$ and \eqref{Tab}, it follows that
 \be
 \big[T\ot{\rm id}_{B(L_2(G))}\big(\pi_{\a}(x)\big)\big]_{ij}
 &=& T\big(\big(\pi_{\a}(x)\big)_{ij}\big)
 =\int_{G}T(\a^{-1}_{h} (x))\bar f_{i}(h)f_{j}(h)dh\\
 &=& \int_{G} \b_{h}^{-1} (Tx)\bar f_{i}(h)f_{j}(h) dh
 =\big(\pi_{\b}(Tx)\big)_{ij}\,.
 \ee
On the other hand,
 $$\l(g) ={\rm id}_{H} \ot \ell(g),$$
where $\ell(g)\, :\, L_2(G) \to L_2(G)$ is the translation by $g.$
Hence, the matrix of $\l(g)$ in $(f_{i})_{i \in I}$ is ${\rm
id}_{H} \ot (a_{ij}),$ where $(a_{ij})_{i,j \in I}$ is a bounded
scalar matrix. Therefore, the matrix of $\l(g) \pi_{\a}(x)$ is
given by
 $$\big(\l(g)\pi_{\a}(x)\big)_{ij}=\sum_{k \in I}
 a_{ik}\big(\pi_{\a}(x)\big)_{kj}.$$
Thus, we deduce that the coefficient at the position $(i, j)$ of
the matrix of $T\ot{\rm id}_{B(L_2(G))}\big(\l(g)\pi_{\a}
(x)\big)$ is
 \be
 \big[T\ot{\rm id}_{B(L_2(G))}\big(\l(g)\pi_{\a}(x)\big)\big]_{ij}
 &=& \sum_{k} a_{ik}\, T \big[(\pi_{\a}(x))_{kj}\big]\\
 &=& \sum_{k}a_{ik}\,\big(\pi_{\b}(T x)\big)_{kj}\\
 &=&\big[\l(g)\pi_{\b}(Tx)\big]_{ij}\,.
 \ee
Hence, it follows that
 \beq\label{Tab1}
 T\ot{\rm id}_{B(L_2(G))}\big(\l(g)\pi_{\a}(x)\big)
 =\l(g)\pi_{\b}(Tx),\quad x \in \M,\ g \in G.
 \eeq
Recall that the family of all finite linear combinations of
$\l(g)\pi_{\a}(x),\ g \in G, \ x \in \M$, is a w*-dense involutive
subalgebra of $\M \rtimes_{\a} G.$ Thus by the normality of $T
\otimes {\rm id}_{B(L_2(G))},$ we deduce our claim.

Now set
 $$\wh T = T \otimes {\rm id}_{B(L_2(G))}
 \big|_{\M\rtimes_{\a}G}\,.$$
Then $\wh T\, :\, \M \rtimes_{\a}G \to \N \rtimes_{\b}G$ is
completely bounded and normal. Moreover, $\wh T$ satisfies
\eqref{def of T hat} by virtue of \eqref{Tab1}. Thus $\wh T$ is
the desired extension of $T.$ The uniqueness of $\wh T$ follows
from the normality and \eqref{def of T hat}.

Let us check the other properties of $\wh T.$ Let $g,\, h \in G,\,
x \in \M.$ Then by \eqref{com pi-l},
 $$\l(g) \pi_{\a}(x)\l(h)
 = \l(g)\l(h) \pi_{\a} (\a_{-h}(x))
 = \l(g+h)\pi_{\a}(\a_{-h} (x)).$$
Thus by \eqref{def of T hat} and \eqref{Tab},
 \be
 \wh T\big(\l(g)\pi_{\a}(x)\l(h)\big)
 &=&\l(g+h)\pi_{\b}\big(T(\a_{-h}(x))\big)\\
 &=&\l(g+h)\pi_{\b}\big(\b_{-h} (Tx)\big)
 = \l(g)\pi_{\b} (T x)\l(h).
 \ee
This yields \eqref{def of T hat1} in the case where $a = \l(g)$
and  $b = \l(h)$ for any $g, h \in G.$ The general case then
follows from the normality of $\wh T.$

(ii) is a consequence of \eqref{dual action} and \eqref{Tab}. If
$T$ is a homomorphism, $\ast$-homomorphism or completely positive
map, then so is $T \otimes {\rm id}_{B(L_2(G))}.$ Thus we get
(iii).

Under the conditions of (iv), $\N \rtimes_{\b}G$ is a subalgebra
of $\M \rtimes_{\a}G.$ If $T$ is a (faithful) conditional
expectation, so is $T \otimes {\rm id}_{B(L_2(G))}.$ Hence $\wh T$
is also a (faithful) conditional expectation.

(v) follows from \eqref{modular gp of dual weight} and
\eqref{Tab}. To prove (vi) we extend, by normality, both $T$ and
$\wh T$ to the extended positive parts of $\M$ and $\M
\rtimes_{\a} G$, respectively. Let $\Phi_\a$ and $\Phi_\b$ be the
operator-valued weights from $\M \rtimes_{\a} G$ to $\M$ and from
$\N \rtimes_{\b} G$ to $\N$, respectively. Then by \eqref{o-v
weight} and the normality we obtain $\Phi_\b\circ\wh T=\wh
T\circ\Phi_\a$. Thus
 \be
 \wh\psi\circ\wh T
 &=&\psi\circ\pi_\b^{-1}\circ\Phi_\b\circ\wh T
 =\psi\circ\pi_\b^{-1}\circ\wh T\circ\Phi_\a\\
 &=&\psi\circ\wh T\circ\pi_\a^{-1}\circ\Phi_\a
 =\psi\circ T\circ\pi_\a^{-1}\circ\Phi_\a\\
 &\le& \f\circ\pi_\a^{-1}\circ\Phi_\a=\wh\f.
 \ee
Hence (vi) is proved. Therefore, the proof of the theorem is
complete. \cqd

\begin{rk}
 It is easy to check that the extension in
Theorem \ref{extension crossed} satisfies the following functorial
property. Let ${\cal L}$ be a third von Neumann algebra, $\gamma$
an automorphic representation of $G$ on ${\mathcal L}$ and $S\,
:\, \N \to {\mathcal L}$ a completely bounded normal map such that
 $S \circ\b_{g} = \gamma_{g} \circ S$ for all $ g \in G$.
  Then $\wh{S \circ T} = \wh S \circ \wh T.$
\end{rk}

We end this section by specializing the extension in Theorem
\ref{extension crossed} to the situation described in section
\ref{Reduction}. Let $\M$ and $\N$ be two von Neumann algebras on
$H$ equipped with two normal faithful states $\f$ and $\psi$,
respectively. Let $G = \bigcup_{n \geq 1}\, 2^{-n}\ent.$ We keep
all the notation in section \ref{Reduction} for both $(\M, \f)$
and $(\N, \psi).$ Set $\R = \M \rtimes_{\s^{\f}} G$ and $\S = \N
\rtimes_{\s^{\psi}} G$. By Theorem \ref{red s-finite}, we have two
increasing sequences $(\R_{n})_{n\ge 1}$ and
 $(\S_{n})_{n\ge 1}$ of finite von Neumann subalgebras of
 $\R$ and $\S$,
respectively, which satisfy all properties there. The
corresponding conditional expectations from  $\R$ onto $\R_{n}$,
respectively, from  $\S$ onto $\S_{n}$, are denoted by $\Phi_{n}$
and $\Psi_n$.

\begin{prop}\label{extension crossed bis}
 Let $T\, :\, \M \to \N$ be a completely bounded
normal map such that
 $$T \circ \s^{\f}_{t}
 = \s_{t}^{\psi} \circ T, \quad t
 \in\real.$$
Let $\wh T\, :\, \R \to \S$ be the extension of $T$ given by
Theorem \ref{extension crossed}. Then
 \begin{enumerate}[\rm(i)]
 \item $\wh T \circ \Phi_{n}
 = \Psi_{n} \circ \wh T$ for every $n$; consequently,
 $\wh T (\R_{n}) \subset \S_{n}$.
 \item Assume in addition that $T\ge0$
and $\psi\circ T\le\f$. Then $\psi_n\circ\wh T\le\f_n$ for every
$n\in\nat$, where $\f_n$ and $\psi_n$ are, respectively, the
states relative to $\f$ and $\psi$ defined by \eqref{def of fn}.
 \end{enumerate}
 \end{prop}

\pf By the definition of $\Phi_n$ in section \ref{Reduction} (see
the proof of Lemma \ref{construction of Rn}) and \eqref{formula od
sfn}, we have
 $$\Phi_n(x)=\int_0^1e^{-ita_n}\s_t^{\wh\f}(x)e^{ita_n}\,dt$$
and a similar formula for $\Psi_n$ with $\f$ replaced by $\psi$.
We then deduce (i) by virtue of \eqref{def of T hat1} and
\eqref{That sfhat}. To prove (ii) recall that $\f_n=e^{-a_n}\wh\f$
and  $\psi_n=e^{-a_n}\wh\psi$. By the fact that $a_{n} \in
\Z(\R_{\wh\psi})$, $a_{n} \in \Z(\R_{\wh\f})$ (see Lemma \ref{bn})
and Theorem \ref{extension crossed}, (i), (vi), we deduce that for
$x\in\R_+$,
 \be
 \psi_n\circ\wh T(x)
 =\wh\psi(e^{-a_n}\wh T(x))
 =\wh\psi(\wh T(e^{-a_n/2}xe^{-a_n/2}))
 \le\wh\f(e^{-a_n/2}xe^{-a_n/2})
 =\f_n(x).
 \ee
This finishes the proof. \cqd


\section{Extensions of maps to noncommutative $L_p$-spaces}
 \label{Extensions of maps to noncommutative Lp-spaces}


In this section we deal with the problem of how to extend a map
between two von Neumann algebras to their noncommutative
$L_p$-spaces. We consider only the $\s$-finite case. Let $\M$ and
$\N$ be two von Neumann algebras equipped with  normal faithful
states $\f$ and $\psi,$ respectively. Let $D_{\f}$ denote the
Radon-Nikod\'ym derivative of the dual weight $\wh\f$ on $\M
\rtimes_{\s^\f}\real$ with respect to the canonical {\it n.s.f.}
trace $\tau_{\f}$ of $\M \rtimes_{\s^\f}\real$. $D_{\psi}$ has the
same meaning relative to $(\N, \psi)$. Consider a positive map
$T\,:\, \M \to \N$  such that for some positive constant $C_1$,
 \beq\label{Linfty bound}
 \psi (T(x)) \leq C_{1}\f(x),\quad
  x \in \M_+.
  \eeq
Given  $1\le p<\8$ define
 $$\begin{array}{llll}
 T_p\, : & D^{1/2p}_{\f}\,\M\, D^{1/2p}_{\f} & \to
 & D^{1/2p}_{\psi}\,\N\, D^{1/2p}_{\psi},\\
 &D^{1/2p}_{\f}x D^{1/2p}_{\f} &\mapsto
 & D^{1/2p}_{\psi}T(x) D^{1/2p}_{\psi}.
 \end{array}$$
Recall that by \cite[Lemma~1.1]{jx-burk}, $D^{1/2p}_{\f}\M
D^{1/2p}_{\f}$ is a dense subspace of $L_p(\M,\f)$. The main
result of this section is the following.

\begin{thm}\label{extension Lp}
 The map $T_{p}$ above extends to a positive bounded map
from $L_p(\M, \f)$ into $L_p(\N, \psi)$ for all $1\le p<\infty$.
Moreover,
 $$\| T_{p} \|\leq C_\8^{1-1/p}
 C_{1}^{1/p}\,,\quad\mbox{where}\  C_\8=\|T(1)\|_{\8}\,.$$
  \end{thm}

Note that there does not exist any additional factor before $C_1$
in the above bound $C_\8^{1-{1/p}}C_{1}^{1/p}$ for the norm of
$T_p$. This is very important for applications, for instance, for
those applications to noncommutative ergodic theory (see section
\ref{Applications to noncommutative maximal inequalities} below).
In the tracial case, Theorem \ref{extension Lp} was proved in
\cite{ye1} with $4C_1$ instead of $C_1$ in the previous estimate
on $\|T_p\|$. Theorem \ref{extension Lp} was announced in
\cite{gold-linKMS}. The proof there presents, unfortunately, a
serious gap.

\smallskip

Now we proceed to the proof of  Theorem \ref{extension Lp}. The
main difficulty lies in the extension of $T_1$. Once this is done,
that of $T_p$ will then follow from a rather easy interpolation
argument via Kosaki's interpolation theorem. For the extension of
$T_1$ we need the following lemma, which is a reformulation of
Lemma 1.2 from \cite{haag-nw} into the present setting. We include
a proof for the convenience of the reader. $\M_h$ denotes the
subspace of selfadjoint elements of $\M$.

\begin{lem}\label{L1 by f}
 Let $x\in \M_h$. Then
 $$ \|D_\f^{1/2}xD_\f^{1/2}\|_1=\inf\big\{\f(a)+\f(b)\;:\; x=a-b,\
 a,\; b\in\M_+\big\}.$$
 \end{lem}

 \pf Denote the infimum on the right-hand side by
 $\rho(x)$. Then $x\mapsto \rho(x)$ defines a seminorm on
 $\M_h$. By \eqref{state vs trace}, for any $x\in\M_+$,
  $$\rho(x)=\f(x)=\tr(D^{1/2}xD^{1/2})=\|D^{1/2}xD^{1/2}\|_1\,.$$
Here and during this proof we denote $D_\f$ simply by $D$.
 Thus it follows that
 $$ \|D^{1/2}xD^{1/2}\|_1\le \rho(x),\quad x\in\M_h\,.$$
 To prove the converse inequality, we  fix $x_0\in\M_h$.
 Then by the Hahn-Banach theorem there exists a
 linear functional $f:\M_h\to\real$ such that
  $$f(x_0)=\rho(x_0)\quad\mbox{and}\quad |f(x)|\le
  \rho(x),\quad\forall\; x\in \M_h\,.$$
We extend $f$ to a complex linear functional on the whole $\M$ by
complexification, still denoted by $f$. Then $f$ is Hermitian and
$-\f\le f\le\f$ on $\M_+$. Thus by the Cauchy-Schwarz inequality,
for $x,\;y\in\M$ we have
   \be
   |f(y^*x)|
   &\le& {\frac12}\big[|(\f+f)(y^*x)| + |(\f-f)(y^*x)|\big]\\
   &\le& \frac12\,\big[\big((\f+f)(y^*y)\big)^{1/2}\big((\f+f)(x^*x)\big)^{1/2}
    +\big((\f-f)(y^*y)\big)^{1/2}\big((\f-f)(x^*x)\big)^{1/2}\big]\\
   &\le& \big(\f(y^*y)\big)^{1/2}\big(\f(x^*x)\big)^{1/2}
   =\|xD^{1/2}\|_2\,\|yD^{1/2}\|_2\,.
   \ee
By the density of $\M D^{1/2}$ in $L_2(\M,\f)$,  we deduce that
there exists a contraction $B$ on $L_2(\M,\f)$ such that
   $$\langle B(xD^{1/2}),\; yD^{1/2}\rangle = f(y^*x),\quad x,
   y\in\M\,.$$
Now  regarding $\M$ as acting standardly on $L_2(\M,\f)$ by  left
multiplication, we claim that $B$ belongs to the commutant of
$\M$. Indeed, for $a\in\M$,
   \be
   \langle Ba(xD^{1/2}),\; yD^{1/2}\rangle
   &=&\langle B(axD^{1/2}),\; yD^{1/2}\rangle=f(y^*ax)\\
   &=&\langle B(xD^{1/2}),\; a^*yD^{1/2}\rangle
   =\langle aB(xD^{1/2}),\; yD^{1/2}\rangle\,.
   \ee
Therefore, $Ba=aB$, so our claim follows. Thus $B$ coincides with
the right multiplication on $L_2(\M,\f)$ by an element $b\in\M$.
Hence we deduce that
    $$f(y^*x)=\langle xD^{1/2}b,\; yD^{1/2}\rangle
    =\tr(D^{1/2}y^* xD^{1/2}b),\quad x,y\in\M\,.$$
Consequently,
    $$\rho(x_0)=f(x_0)=\tr(D^{1/2}x_0D^{1/2}b)
    \le \|D^{1/2}x_0D^{1/2}\|_1 \|b\|_\8\le
    \|D^{1/2}x_0D^{1/2}\|_1\,.$$
Therefore, the lemma is proved.\cqd

\begin{lem}\label{L1 extension}
 $T_{1}$ extends to a positive bounded map from $L_1(\M,
\f)$ into $L_1(\N, \psi)$ with norm majorized by $C_{1}.$
 \end{lem}

\pf  Let $x\in\M$ and $y = D_{\f}^{1/2}xD_{\f}^{1/2}$. Assume
first $x \geq 0$. Then $y \geq 0$, so $ T_{1} (y) \geq 0 $ for $T$
is positive. Hence, by \eqref{state vs trace} and \eqref{Linfty
bound},
 \be
 \|T_{1}(y)\|_{L_1(\N, \psi)}
 &=&\tr(T_{1}(y))=\tr(D_{\psi}^{1/2}T(x)D_{\psi}^{1/2})\\
 &=&\psi (T(x))\leq C_{1}\f(x)
 =C_{1}\|y\|_{L_1(\M, \f)}.
 \ee
Now assume that $x$ is selfadjoint  and $\e>0$.  Then by Lemma
\ref{L1 by f}, there exist $a, b\in\M_+$ such that
 $x=a-b$ and
  $$\|D_{\f}^{1/2}aD_{\f}^{1/2}\|_{L_1(\M,\f)}
 +\|D_{\f}^{1/2}bD_{\f}^{1/2}\|_{L_1(\M,\f)}
 \le \|y\|_{L_1(\M,\f)}+\e.$$
It follows that
 $$T_{1}(y)=D_{\psi}^{1/2}T(a)D_{\psi}^{1/2}
 -D_{\psi}^{1/2}T(b)D_{\psi}^{1/2}$$
and
 \be
 \|T_{1}(y)\|_{L_1(\N, \psi)}
 &\le& \|D_{\psi}^{1/2}T(a)D_{\psi}^{1/2}\|_{L_1(\N, \psi)}
  + \|D_{\psi}^{1/2}T(b)D_{\psi}^{1/2}\|_{L_1(\N,
  \psi)}\\
 &\le& C_1\big(\|D_{\f}^{1/2}aD_{\f}^{1/2}\|_{L_1(\M, \f)}
  + \|D_{\f}^{1/2}bD_{\f}^{1/2}\|_{L_1(\M,
  \f)}\big)\\
 &\le& C_1\big(\|y\|_{L_1(\M, \f)}+\e\big),
 \ee
whence again
  $$\|T_{1}(y)\|_{L_1(\N, \psi)}
  \le C_1\|y\|_{L_1(\M,\f)}\,.$$
Finally, decomposing any $x \in \M$ into its real and imaginary
parts, we get
 $$\|T_{1} (y) \|_{L_1(\N,\psi)}
 \leq 2\, C_{1} \| y\|_{L_1(\M,\f)}.$$
Therefore, $T_{1}$ is bounded relative to the $L_1$-norms. Since
$D_{\f}^{1/2}\M D_{\f}^{1/2}$ is dense in $L_1(\M, \f),$ $T_{1}$
extends to a bounded map from $L_1(\M, \f)$ into $L_1(\N, \psi)$
with $\|T_1\| \leq 2C_1.$ Since $T_1$ is positive, so is its
extension (which is denoted again by $T_1$).

Thus it remains to drop the factor 2 from the previous estimate on
$\|T_1\|$. To this end, we consider the adjoint: $T_1^*: \N\to\M$.
Since $T_1^*$ is positive, $T_1^*$ attains its norm at the
identity of $\N$ (see \cite{pau-cb}). Therefore,
 $$\|T_1\|=\|T_1^*\|=\|T_1^*(1)\|_\8\,.$$
Hence, we are reduced to showing $\|T_1^*(1)\|_\8\le C_1$. This is
easy by duality. Indeed, let $x\in\M_+$ and
$y=D_\f^{1/2}xD_\f^{1/2}$. Then by \eqref{state vs trace} and
\eqref{Linfty bound},
 $$\tr(T_1^*(1)\, y)=\tr(T_1(y))
 =\tr\big(D_\psi^{1/2}T(x)D_\psi^{1/2}\big)
 =\psi(T(x))\le C_1\f(x)=C_1\|y\|_{L_1(\M,\f)}\,.$$
Since $D_\f^{1/2}\M_+D_\f^{1/2}$ is dense in the positive cone of
$L_1(\M,\f)$, we deduce the desired estimate on $\|T_1^*(1)\|_\8$.
Hence the proof of the lemma is complete. \cqd

\medskip

\noindent {\it Proof of Theorem \ref{extension Lp}.} To prove that
$T_{p}$ extends to a bounded map from $L_p(\M, \f)$ to $L_p(\N,
\psi)$ for all $1 < p < \infty,$ we will use interpolation. We
consider the following symmetric injection of $\M$ into $L_1(\M,
\f)$:
 $$\begin{array}{llll}
 j_{\f}\, : & \M & \to & L_1(\M, \f),\\
 &x &\mapsto & D^{1/2}_{\f}x D^{1/2}_{\f},
 \end{array}$$
Then $j_{\f}$ turns $(\M,\; L_1(\M,\f))$ into a compatible couple,
so we can consider their complex interpolation space $(\M,\;
L_1(\M, \f))_{1/p}$ for any $1<p<\8$. By \cite[Theorem
9.1]{kos-int}, this space can be isometrically identified with
$L_p(\M,\f)$. More precisely,  define
 $j^{p}_{\f} (x)
 =D_{\f}^{1/2p}xD_{\f}^{1/2p},\ x \in \M.$ Then
 $j^{p}_{\f}$ extends to an isometry from
$(\M,\; L_1(\M, \f))_{1/p}$ onto $L_p(\M,\f)$. For $(\N, \psi)$ we
use similar notation. Under the injections $j_{\f}$ and
$j_{\psi},$  $T_1$ is viewed as the same map as $T$ (on the
intersection space $\M$). Therefore, by Lemma \ref{L1 extension}
and interpolation, $T$ is bounded from $(\M,\; L_1(\M, \f))_{1/p}$
into $(\N,\; L_1(\N, \psi))_{1/p}$ with norm majorized by
$C_{\8}^{1 - {1/p}} C_{1}^{1/p}$ (recalling that
$C_\8=\|T\|_{\M\to\N}$). From this, and using the isometries
$j^{p}_{\f},\ j^{p}_{\psi}$ defined above, we deduce that for any
$x \in \M$,
 \be
 \big\|T_p\big(D_{\f}^{1/2p}xD_{\f}^{1/2p}\big)
 \big\|_{L_p(\N, \psi)}
 &=&\big\|j_{\psi}^p(T(x))\big\|_{L_p(\N, \psi)}
 =\big\|T(x)\big\|_{(\N,\; L_1(\N, \psi))_{1/p}}\\
 &\le& C_{\8}^{1 - 1/p}\,C_{1}^{1/p}
 \big\|x\big\|_{(\M,\; L_1(\M, \f))_{1/p}}\\
 &=&C_{\8}^{1 - 1/p}\,C_{1}^{1/p}
 \big\|j_{\f}^p(x)\big\|_{L_p(\M, \f)}\\
 &=&C_{\8}^{1 - 1/p}\,C_{1}^{1/p}
 \big\|D_{\f}^{1/2p}xD_{\f}^{1/2p}\big\|_{L_p(\M, \f)}\,.
 \ee
Therefore, by the density of $D_{\f}^{1/2p}\M\, D_{\f}^{1/2p}$ in
$L_p(\M, \f)$, this implies that $T_p$ extends to a bounded map
from $L_p(\M,\f)$ to $L_p(\N,\psi)$ with norm controlled by
$C_{\8}^{1 - 1/p}\,C_{1}^{1/p}$. Since $T_p$ is positive on
$D_{\f}^{1/2p}\,\M\, D_{\f}^{1/2p}$, so is its extension on
$L_p(\M,\f)$. \cqd

\smallskip

Let $T$ be as in Theorem \ref{extension Lp}. The extension of
$T_p$ will again be denoted by $T_p\,$. Consider the adjoint map
of $T_1$: $S=T_1^*\,:\, \N\to\M$.  We claim that $S$ satisfies the
same conditions as $T$. Namely, $S$ is positive and
 \beq\label{S Linfty}
 \f(S(y))\le C_\8 \psi(y),\quad y\in\N_+\,.
 \eeq
Indeed, the positivity of $S$ was already observed during the
proof of Lemma \ref{L1 extension}.  On the other hand, for
$y\in\N_+$ we have
  \be
  \f(S(y))
  &=&\tr(S(y)D_{\f})
  =\langle T_1^*(y),\;D_{\f}\rangle
  =\langle y,\;D_{\psi}^{1/2}T(1)D_{\psi}^{1/2}\rangle\\
  &=&\tr\big(D_{\psi}^{1/2}yD_{\psi}^{1/2}\,T(1) \big)
  \le\|T(1)\|_\8\psi(y)=C_\8\,\psi(y).
  \ee
Therefore, applying Theorem \ref{extension Lp} to $S$, we get the
extension $S_p\,:\, L_p(\N,\psi)\to L_p(\M,\f)$. It is easy to
check that $S_1^*=T$ and $S_p^*=T_{p'}$  for any $1<p<\8$ ($p'$
being conjugate to $p$). Consequently, $T$ is normal.

\smallskip

We record the discussion above in the following.

\begin{prop}\label{adjoint of T}
 Let $T$ and $T_p$ be as in Theorem \ref{extension Lp}
$(T_p$ also denoting the extension$)$. Let $S=T_1^*$.
 \begin{enumerate}[\rm(i)]
 \item The  map $S\,:\, \N\to\M$ is
positive and satisfies \eqref{S Linfty}.
 \item Let $S_p: L_p(\N,\psi)\to
L_p(\M,\f)$ be the extension of $S$. Then $S_1^*=T$ and
$S_p^*=T_{p'}$ for every $1< p<\8$.
 \item  $T$ is normal.
 \end{enumerate}
 \end{prop}

The extension in Theorem \ref{extension Lp} is  symmetric. We
could also consider the left extension: $xD_\f^{1/p}\mapsto
T(x)D_\psi^{1/p}$ ($x\in\M$). More generally, for any
$0\le\theta\le1$ we can define
 $$\begin{array}{llll}
 T_{p,\th}\,: & D^{(1-\th)/p}_{\f}\,\M\,
 D^{\th/p}_{\f} & \to
 & D^{(1-\th)/p}_{\psi}\,\N\, D^{\th/p}_{\psi},\\
 &D^{(1-\th)/p}_{\f}x D^{\th/p}_{\f} &\mapsto
 & D^{(1-\th)/p}_{\psi}T(x) D^{\th/p}_{\psi}.
 \end{array}$$
Thus $T_{p,1/2}$ agrees with $T_p$. We do not know whether
$T_{p,\th}$ extends to a bounded map from $L_p(\M,\f)$ to
$L_p(\N,\psi)$. However, if in addition $T$ satisfies
 \beq\label{T com sf}
 T \circ \s^{\f}_{t} = \s_{t}^{\psi} \circ T, \quad t
 \in\real,
 \eeq
then  $T_{p,\th}$ indeed extends.  Recall that $\M_a$ (resp.
$\N_a$) denotes the family of all analytic elements of $\M$ with
respect to $\s^\f$ (resp. $\N$ relative to $\s^\psi$). By
\cite[Lemma~1.1]{jx-burk}, $D^{(1-\th)/p}_{\f}\,\M_a\,
 D^{\th/p}_{\f}$ is dense in $L_p(\M,\f)$.

\begin{prop}\label{extension Lp bis}
 Let $T\,:\, \M\to\N$ satisfy \eqref{T com sf}. Then
 $T_{p,\theta}=T_p$
on $D^{(1-\th)/p}_{\f}\,\M_a\,
 D^{\th/p}_{\f}$. Consequently, if in addition, $T$ is positive
and satisfies \eqref{Linfty bound}, then $T_{p,\theta}$ extends to
a bounded map from $L_p(\M, \f)$ to $L_p(\N, \psi)$ and its
extension coincides with that of $T_p$ in Theorem \ref{extension
Lp}.
 \end{prop}

\pf We first note that by \eqref{T com sf}, $T$ maps $\M_{a}$ into
$\N_{a}$. On the other hand, it is easy to see that
 $$D^{(1-\th)/p}_{\f}\,\M_{a}\,
 D^{\th/p}_{\f}=D^{1/2p}_{\f}\,\M_{a}\,
 D^{1/2p}_{\f}.$$
Now let $x \in \M_a$. Then
 $$
 D_{\f}^{(1-\th)/p}xD_{\f}^{\th/p}
 =\s^{\f}_{-\frac{i(1-\theta)}p}(x)D_{\f}^{1/p}\, .$$
Therefore, by \eqref{T com sf}, we get
 \be
 T_{p,1}\big(\s^{\f}_{-\frac{i(1-\theta)}{p}}(x)D_{\f}^{1/p}\big)
 &=&T\big(\s^{\f}_{-\frac{i(1-\theta)}{p}}(x)\big)\,D_{\psi}^{1/p}
 = \s^{\psi}_{-\frac{i(1-\theta)}{p}}\big(T(x)\big)D_{\psi}^{1/p}\\
 &=&D_{\psi}^{(1-\th)/p}T(x)D_{\psi}^{\th/p}
 =T_{p, \th}\big(D_{\f}^{(1-\th)/p}xD_{\f}^{\th/p}\big)\,.
 \ee
This proves the first part of the proposition. The  second simply
follows from Theorem \ref{extension Lp}. \cqd

\begin{rk}
 Theorem \ref{extension Lp} can be extended to the weighted case too.
Let $\M$ and $\N$ be two general von Neumann algebras equipped
with {\it n.s.f.} weights $\f$ and $\psi,$ respectively. Let
$T\,:\, \M \to \N$ be  a positive map such that for some positive
constant $C_1$,
 \be
 \psi (T(x)) \leq C_{1}\f(x),\quad
  x \in \M_+.
  \ee
Consider
 $$\begin{array}{llll}
 T_p\, : & D^{1/2p}_{\f}\,\rM_\f\, D^{1/2p}_{\f} & \to
 & D^{1/2p}_{\psi}\,\rM_\psi\, D^{1/2p}_{\psi},\\
 &D^{1/2p}_{\f}x D^{1/2p}_{\f} &\mapsto
 & D^{1/2p}_{\psi}T(x) D^{1/2p}_{\psi}.
 \end{array}$$
Then $T_{p}$  extends to a  bounded map from $L_p(\M, \f)$ into
$L_p(\N, \psi)$ for all $1\le p<\infty$.
 \end{rk}

The proof of this statement is essentially the same as that of
Theorem \ref{extension Lp}. Now instead of Kosaki's interpolation
theorem, we use that of Terp \cite{terp-int}. Note also that
$D^{1/2p}_{\f}\,\rM_\f\, D^{1/2p}_{\f}$ is dense in $L_p(\M, \f)$
for any $1\le p<\8$ (see  \cite{gold-linKMSS} and
\cite{terp-int}).

\smallskip

\n{\bf Convention.} In the sequel, we will denote, by the same
symbol $T$, all maps $T_{p}$ and $T_{p,\theta}$ as well as their
extensions between the $L_p$-spaces in Theorem \ref{extension Lp}
and Proposition \ref{extension Lp bis}, whenever no confusion can
occur. This is consistent with Kosaki's interpolation theorem. See
the discussion in the proof of Theorem \ref{extension Lp}.

\medskip

We end this section by three examples, which are special cases of
Theorem \ref{extension Lp}. In all three, the map $T$ satisfies
both \eqref{Linfty bound} and \eqref{T com sf}, so the extensions
to the $L_p$-spaces can be made from $T_p$ in Theorem
\ref{extension Lp} or any $T_{p,\theta}$ in Proposition
\ref{extension Lp bis}.

\begin{ex}\label{extension Lp tracial}
 The first example is the tracial case, i.e., when both $\f$ and
$\psi$ are tracial.
 Then it is well known that any positive map $T\,:\,
\M\to\N$ satisfying (5.1) extends to a bounded map between the
usual noncommutative $L_p$-spaces constructed from the traces $\f$
and $\psi$. This is very easy to prove by interpolation (cf.
\cite{ye1}). Note that in \cite{ye1} the estimate on $\|T_1\|$ is
$4C_1$. Also note that in this case, \eqref{T com sf} is trivially
satisfied.
 \end{ex}

\begin{ex}\label{extension Lp cond}
The second example  concerns conditional expectations. Let $\N$ be
a von Neumann subalgebra of $\M$ and $\psi=\f\big|_\N\,$. Let
$\E\,:\, \M\to\N$ be a normal faithful conditional expectation
such that $\f\circ\E=\f$. Then it is well known that $\E$ commutes
with the modular automorphism group $\s^{\f}_t$ (cf.
\cite[1.4.3]{con-classi}). Thus $\E$ extends to a bounded map
between the noncommutative $L_p$-spaces for all $1\le p<\8$. In
fact, in this special case, Theorem \ref{extension Lp} becomes
\cite[Lemma~2.2]{jx-burk}. Note that the extension of $\E$ between
the $L_p$-spaces possesses all the usual properties of a
conditional expectation as in the commutative case. In particular,
it is an $\N$-bimodular contractive projection from $L_p(\M,\f)$
onto $L_p(\N,\psi)$.
 \end{ex}

\begin{ex}\label{extension Lp isomorphism}
 The third example is that when $T$ is an isomorphism preserving
states, i.e., $\psi\circ T=\f$. In this case condition \eqref{T
com sf} is again automatically satisfied. Indeed, it is easy to
check that $T^{-1}\circ\s_t^\psi\circ T$ is an automorphism group
satisfying the KMS condition relative to $\f$, so it coincides
with $\s_t^\f$. Also in this case Lemma \ref{L1 extension} above
admits a straightforward proof as follows. Let $x\in\M$ and
$y\in\N_a$. Then
 \be
 \tr\big(y\,D_\psi^{1/2}T(x)D_\psi^{1/2}\big)
 &=&\tr\big(D_\psi^{1/2}\s_{i/2}^\psi(y)T(x)D_\psi^{1/2}\big)
 =\psi\big(\s_{i/2}^\psi(y)T(x)\big)\\
 &=&\f\big(\s_{i/2}^\f(T^{-1}(y))\,x\big)
 =\tr\big(T^{-1}(y)\,D_\f^{1/2}xD_\f^{1/2}\big);
 \ee
so
 $$\big|\tr\big(y\,D_\psi^{1/2}T(x)D_\psi^{1/2}\big)\big|
 \le \|T^{-1}(y)\|_\8\,\|D_\f^{1/2}xD_\f^{1/2}\|_1
 \le \|y\|_\8\,\|D_\f^{1/2}xD_\f^{1/2}\|_1.$$
This implies the boundedness of  $T_1$ on $D_\f^{1/2}\M
D_\f^{1/2}$.
 \end{ex}


\section{Applications to noncommutative martingale
inequalities}
  \label{Applications to noncommutative martingale
theory}


Since the establishment of the noncommutative Burkholder-Gundy
inequalities in \cite{px-BG}, the theory of noncommutative
martingale inequalities has been rapidly developed. Many of the
classical inequalities in the usual martingale theory have been
transferred into the noncommutative setting. We refer, for
instance,  to \cite{ju-doob} for the  Doob maximal inequality, to
\cite{jx-burk, jx-ros} for the Burkholder/Rosenthal inequalities,
to \cite{ran-mtrans, ran-weak, ran-cond} for several weak type
$(1, 1)$ inequalities and to \cite{par-ran-gu} for the Gundy
decomposition.

The objective of this section and the next one is to show how to
use Theorem \ref{red s-finite} to reduce all these inequalities in
the nontracial case to those in the tracial one. As a consequence,
their best constants in the general case coincide with the
corresponding ones in the tracial case. This section deals with
square function type inequalities. The Doob maximal inequality,
the only exception among those quoted previously, is postponed to
the next section, where the maximal ergodic inequalities will also
be considered.


\subsection{The framework}
 \label{The framework}

Throughout this and the next sections, $\M$ will denote a von
Neumann algebra equipped with a distinguished normal faithful
state $\f$ and $\s=\s^\f$ the modular automorphism group of $\f$.
We denote $L_p(\M,\f)$ simply by $L_p(\M)$. Assume that $\N$ is a
von Neumann subalgebra of $\M$ and that there exists a normal
faithful conditional expectation $\E$ from $\M$ onto $\N$ such
that
 \beq\label{f-preserving cond}
 \f\circ\E=\f\,.
 \eeq
It is well known that such an $\E$ is unique and satisfies
 \beq\label{mod gp com cond}
 \s_t\circ\E=\E\circ\s_t\,,\quad t\in\real.
 \eeq
Recall that the existence of  a conditional expectation $\E$
satisfying \eqref{f-preserving cond} is equivalent to the
$\s^\f$-invariance of $\N$ (see \cite{tak-cond}). The restriction
of $\s^\f$ to $\N$ is the modular automorphism group of
$\f\big|_\N$. We will not distinguish $\f$ and $\s^\f$ from  their
restrictions to $\N$.

\smallskip

Let $G=\bigcup_{n\ge 1} 2^{-n}\ent$ (fixed throughout this and the
next sections).  For notational simplicity set
 $$\R(\M) = \M \rtimes_{\s} G\quad\mbox{and}\quad
 \R(\N) = \N\rtimes_{\s} G\,.$$
By virtue of \eqref{mod gp com cond}, $\R(\N)$ is naturally viewed
as a von Neumann subalgebra of $\R(\M)$. Let $\wh \f$ be the dual
weight on $\R(\M)$ of $\f$ (recalling that $\wh\f$ is a normal
faithful state). Its restriction to $\R(\N)$ is the dual weight of
$\f \big|_{\N}$. As usual, we denote this restriction again by
$\wh \f$. The increasing sequence of  the von Neumann subalgebras
of $\R(\M)$ constructed in Theorem \ref{red s-finite} relative to
$\M$ is denoted by $(\R_m(\M))_{m\ge 1}$, and that relative to
$\N$ by $(\R_m(\N))_{m\ge 1}$. All $\R_m(\M)$ and $\R_m(\N)$ are
von Neumann subalgebras of $\R(\M)$. From the proof of Theorem
\ref{red s-finite} we easily see that
 $$\R_m(\N)=\R_m(\N)\cap \R_m(\M)\,,\quad m\in\nat.$$

Let $\Phi\,:\, \R(\M)\to \M$ be the conditional expectation
defined by \eqref{cond expect}. Its restriction to $\R(\N)$ is the
corresponding conditional expectation from $\R(\N)$ onto $\N$,
again denoted by $\Phi$. Let $\Phi_m\,:\, \R(\M)\to \R_m(\M)$  be
the conditional expectation constructed in Theorem \ref{red
s-finite}. On the other hand, by virtue of \eqref{mod gp com cond}
and Theorem \ref{extension crossed}, $\E$ extends to a normal
faithful conditional expectation $\wh\E$ from $\R(\M)$ onto
$\R(\N)$. Then by \eqref{dual state}, \eqref{def of T hat} and
\eqref{That sfhat}, we have
 \beq\label{cond hat f hat}
 \wh\f\circ\wh\E = \wh\f
 \quad\mbox{and} \quad
 \s_{t}^{\wh\f}\circ\wh\E
 =\wh\E\circ\s_{t}^{\wh\f}\,.
 \eeq
 By \eqref{cond expect}, \eqref{def of T hat} and
Proposition \ref{extension crossed bis}, we find
 \beq\label{cond hat cond hat}
 \wh\E\circ\Phi=\Phi\circ\wh\E\quad\mbox{and}\quad
 \wh\E\circ\Phi_m=\Phi_m\circ\wh\E\,,\quad m\in\nat.
 \eeq
It follows that $\R(\N)$ and $\R_m(\M)$ are respectively invariant
under $\Phi_m$ and $\wh\E$ for all $m\in\nat$; moreover,
$\Phi_m\big|_{\R(\N)}\,:\, \R(\N)\to \R_m(\N)$ and
$\wh\E\big|_{\R_m(\M)}\,:\, \R_m(\M)\to \R_m(\N)$ are normal
faithful conditional expectations. Again, we will not distinguish
these conditional expectations and their respective restrictions.
By Theorem \ref{red s-finite} and  the discussion above, all
conditional expectations $\Phi,\, \Phi_m$ and $\wh\E$ commute, and
further commute with  $\s^{\wh\f}$; moreover, all these maps
preserve the dual state $\wh\f$. This commutation is shown in the
following diagram:
 \[
 \xymatrix{
 & \M\; \ar[ld]_{\E}\ar[d]_{\s_t^\f} & \R(\M) \ar[l]_{\Phi} \ar[d]_{\s_t^{\wh\f}}
 \ar[r]^{\Phi_m}\; & \R_m(\M)\; \ar[d]_{\s_t^{\wh\f}} \ar[rd]^{\wh\E} \\
 \N\; \ar[rd]_{\s_t^\f}  \;& \M \ar[d]_{\E} & \R(\M)\;
\ar[l]_{\Phi} \ar[d]_{\wh\E}
 \ar[r]^{\Phi_m}\; & \R_m(\M) \ar[d]_{\wh\E}
 & \R_m(\N)\ar[ld]^{\s_t^{\wh\f}} \\
  & \N  & \R(\N)\; \ar[l]_{\Phi} \ar[r]^{\Phi_m}\; & \R_m(\N)
 }\]

\medskip

By Remark~\ref{subLp}, $L_p(\N,\f)=L_p(\N)$ is naturally
identified as a subspace of $L_p(\M)$. On the other hand,
\cite[Lemma~2.2]{jx-burk} (or Example~\ref{extension Lp cond})
implies that $\E$ extends to a positive contractive projection
from $L_p(\M)$ onto $L_p(\N)$ ($1\le p\le\8$), again denoted by
$\E$, which possesses the following modular property. Let $1\le p,
q, r\le\8$ such that
 $1/p+1/q+1/r\le1$. Then
  \beq\label{mod of cond}
  \E(axb)=a\E(x)b\,,\quad a\in L_q(\N),\ b\in L_r(\N),\
  x\in L_p(\M)\,.
  \eeq
In what follows all spaces $L_p(\R(\M)),\ L_p(\R_m(\M)),\
L_p(\R(\N))$ and $L_p(\R_m(\N))$ are relative to $\wh\f$. They are
naturally identified as subspaces of $L_p(\R(\M))$. On the other
hand, since $\wh\f\big|_{\M}=\f$ (see \eqref{dual state}) and
$\M$, $\N$ are $\s^{\wh\f}$-invariant (see \eqref{modular gp of
dual weight}), $L_p(\M)$ and $L_p(\N)$ are subspaces of
$L_p(\R(\M))$ too. Thus, $L_p(\R(\M))$ is the largest space among
all these noncommutative $L_p$-spaces. Since $\Phi,\;\Phi_m$ and
$\wh\E$ preserve $\wh\f$ and commute with $\s^{\wh\f}$,  these
conditional expectations extend to positive contractive
projections on $L_p(\R(\M))$ for all $1\le p\le\8$; moreover,
their extensions satisfy the  modular property \eqref{mod of
cond}. Since all these conditional expectations commute, so do
their extensions. As usual, we use the same symbol to denote a map
and its extensions.

\smallskip
Now we fix  an increasing filtration $(\M_n)_{n\ge 1}$ of von
Neumann subalgebras of $\M$ whose union is w*-dense in $\M$.
Assume that for each $n\ge 1$ there exists a normal faithful
conditional expectation $\E_n$ from $\M$ onto $\M_n$ such that
$\f\circ\E_n=\f$. Then
 \beq\label{filtration of cond}
 \E_n\circ\E_m=\E_m\circ\E_n=\E_{\min(n,m)}\,,
 \quad m,\,n\in \nat.
 \eeq
The preceding discussion applies, of course,  to each $\M_n$ in
place of $\N$. Thus we have the crossed product $\R(\M_n)$ and the
subalgebras $\R_m(\M_n)$. Also, each $\E_n$ extends to a
conditional expectation $\wh\E_n$ from $\R(\M)$ to $\R(\M_n)$. By
\eqref{filtration of cond} and \eqref{def of T hat}, we find
 \beq\label{filtration of cond hat}
 \wh\E_n\circ\wh\E_m=\wh\E_m\circ\wh\E_n=\wh\E_{\min(n,m)}\,,
 \quad m,\ n\in \nat.
 \eeq
All previous assumptions and notation will be kept fixed
throughout this section.


\subsection{Martingale inequalities}


Let $\M$ and $(\M_n)$ be fixed as in the previous subsection. By
definition, an $L_p$-martingale with respect to $(\M_n)_n$ is a
sequence $x=(x_n)\subset L_p(\M)$ ($1\le p\le\8$) such that
  $$\E_n(x_{n+1})=x_n,\quad \forall\; n\in\nat .$$
In this case, $x$ is  adapted in the sense that $x_n\in L_p(\M_n)$
for all $n$.  Define
  $$\|x\|_p=\sup_n\|x_n\|_p\,.$$
If $\|x\|_p<\8$, $x$ is called a bounded $L_p$-martingale. The
martingale difference sequence of $x$ is defined to be
$dx=(dx_n)_{n\ge 1}$ with $dx_n=x_n-x_{n-1}$ ($x_0=0$ by
convention).

\begin{rk}\label{mart conv}
 It is an easy exercice to check the following two
 properties:
 \begin{enumerate}[\rm(i)]
 \item Let $x_\8\in L_p(\M)$ with $1\le p\le\8$, and
let $x_n=\E_n(x_\8)$. Then $x=(x_n)$ is a bounded $L_p$-martingale
and $x_n$ converges to $x_\8$ in $L_p(\M)$ (in the w*-topology for
$p=\8$). Moreover, $\|x\|_p=\|x_\8\|_p$.
 \item Conversely, let $x=(x_n)$ be a bounded $L_p$-martingale
with $1<p\le\8$. Then there exists $x_\8\in L_p(\M)$ such that
$x_n=\E_n(x_\8)$ for all $n$.
 \end{enumerate}
This remark allows us to not distinguish a martingale $x$ and its
final value $x_\8$ whenever the latter exists. This also explains
why we use the letter $x$ to denote sometimes an operator in
$L_p(\M)$, sometimes a martingale. We will also identify a
martingale with its difference sequence. In the sequel all
martingales are with respect to $(\M_n)$ unless explicitly stated
otherwise.
 \end{rk}

Now we can begin to state the noncommutative martingale
inequalities we are interested in. In the sequel, the letters
$\a_{p}, \b_{p}$, $\ldots$ will denote positive constants
depending only on $p$, and $C$ an absolute positive constant. The
simplest noncommutative martingale inequalities are the
noncommutative Khintchine inequalities, which are of paramount
importance in noncommutative analysis.

\medskip\n{\bf Khintchine inequalities.}
Let $(\e_n)$ be a Rademacher sequence  on a probability space
$(\O,\; P)$, i.e., an independent sequence with
$P(\e_n=1)=P(\e_n=-1)=1/2$ for all $n$. Recall that the classical
Khintchine inequality asserts that for any $p<\8$,
 $$\big\|\sum_n \e_na_n\big\|_p\sim
 \big(\sum_n|a_n|^2\big)^{1/2}$$
holds for all finite sequences $(a_n)\subset\com$, where the
equivalence constants depend only on $p$. Using the Fubini theorem
we then deduce that for any finite sequence $(a_n)$ in a
commutative $L_p$-space we have
 $$\big(\mathbb E\big\|\sum_n \e_na_n\big\|_p^p\big)^{1/p}\sim
 \big\|\big(\sum_n|a_n|^2\big)^{1/2}\big\|_p\,,$$
where $\mathbb E$ is the expectation on $\O$.

The noncommutative analogue of the previous equivalence takes,
unfortunately,  a much less simple form due to the existence of
two different absolute values of operators because of the
noncommutativity. We now have two square functions:
 $$\big(\sum_na_n^*a_n\big)^{1/2}\quad \mbox{and}\quad
 \big(\sum_na_na_n^*\big)^{1/2}$$
for any finite sequence $(a_n)\subset L_p(\M)$. Accordingly, we
introduce  the column and row $L_p$-spaces. The column space
$L_p(\M; \el_2^c)$ is the Banach space of all sequences
$a=(a_n)_{n\ge 1}\subset L_p(\M)$ such that
 $$\|a\|_{L_p(\M; \el_2^c)}=
 \big\|\big(\sum_{n\ge1}a_n^*a_n\big)^{1/2}\big\|_{p}<\8.$$
The row space $L_p(\M; \el_2^r)$ consists of all $a$ such that
$a^*\in L_p(\M; \el_2^c)$ and is equipped with the norm
$\|a\|_{L_p(\M;\el_2^r)}=\|a^*\|_{L_p(\M;\el_2^c)}\,$. $L_p(\M;
\el_2^c)$ and $L_p(\M; \el_2^r)$ can be respectively regarded as
the column and row subspaces of $L_p(\M\bar\ot B(\el_2))$ (see
\cite{px-BG} for more details). Then they are $1$-complemented in
$L_p(\M\bar\ot B(\el_2))$ for all $1\le p\le\8$.

Now  define the space $CR_p[L_p(\M)]$ as follows. If $p\ge2$,
 $$CR_p[L_p(\M)]=L_p(\M;\el_2^c)\cap L_p(\M;\el_2^r)$$
equipped with the intersection norm:
 $$\big\|(a_n)\big\|_{CR_p[L_p(\M)]}=
 \max\big\{\big\|(a_n)\big\|_{L_p(\M;\el_2^c)},\;
 \big\|(a_n)\big\|_{L_p(\M;\el_2^r)}\big\}.$$
If $p<2$,
 $$CR_p[L_p(\M)]=L_p(\M;\el_2^c)\,+\, L_p(\M;\el_2^r)$$
equipped with the sum norm:
 $$\big\|(a_n)\big\|_{CR_p[L_p(\M)]}=
 \inf\big\{\big\|(b_n)\big\|_{L_p(\M;\el_2^c)}\,
 +\,\big\|(c_n)\big\|_{L_p(\M;\el_2^r)}\big\},$$
where the infimum runs over all decompositions $a_n=b_n + c_n$
with $b_n, c_n\in L_p(\M)$.

\smallskip

We are now ready to state the noncommutative Khintchine
inequalities. Recall that the vector-valued $L_p$-space $L_p(\O;
L_p(\M))$ can be identified with $L_p(L_\8(\O)\bar\ot\M)$.

\begin{thm}
 Let $1\le p<\8$ and $(x_n)$ be a finite sequence in
$L_p(\M)$. Then
 $$\big\|\sum_n\e_na_n\big\|_{L_p(\O;L_p(\M))}
 \sim\big\|(a_n)\big\|_{CR_p[L_p(\M)]}\,.\leqno({\rm K}_p)$$
More precisely, there exist two universal positive constants $A$
and $B$ such that
 \be
 A^{-1}\,\big\|(a_n)\big\|_{CR_p[L_p(\M)]}\le
 \big\|\sum_n\e_na_n\big\|_{L_p(\O;L_p(\M))}
 \le\big\|(a_n)\big\|_{CR_p[L_p(\M)]}
 \ee
for $1\le p\le2$ and
 \be
 \big\|(a_n)\big\|_{CR_p[L_p(\M)]}\le
 \big\|\sum_n\e_na_n\big\|_{L_p(\O;L_p(\M))}
 \le B\,\sqrt p\,\big\|(a_n)\big\|_{CR_p[L_p(\M)]}
 \ee
for $2\le p<\8$.
 \end{thm}

\smallskip

The previous inequalities were first proved for $1<p<\8$  and
Schatten classes (i.e., for $\M=B(\el_2)$) in \cite{lust-khin} and
then for $1\le p<\8$ and tracial noncommutative $L_p$-spaces in
\cite{LPP} (see \cite{pis-ast} for the optimal order ${\rm
O}(\sqrt p)$). The arguments of \cite{LPP} work for the type III
case too. See \cite{haag-mu-kh} for a very simple proof with
better constants in the case $p=1$. This last work also provides
the best constants for certain other random variables instead of
$(\e_n)$, including complex Gaussians and type III Fermions. We
refer the interested reader to \cite{ju-araki}, \cite{jx-ros},
\cite{pis-ast} for more Khintchine type inequalities in the
noncommutative setting.

Using Theorem \ref{red s-finite}, one can easily reduce the
general case of $({\rm K}_p)$ to the tracial one. See the proof of
Theorem \ref{red mart} below, notably the part concerning $({\rm
S}_p)$.

\medskip\n{\bf Martingale transforms.} We now
consider the unconditionality of noncommutative martingale
difference sequences, i.e., the noncommutative martingale
transforms by sequences of signs.

\begin{thm}
 Let $1 < p < \infty$. Then for any finite $L_p$-martingale
$x=(x_n)$,
 $$
 \big\| \sum_{n\ge1} \e_{n} dx_{n}\big\|_p
 \le \kappa_{p}\, \| x \|_{p},\quad \forall\;\e_{n}
 = \pm 1. \leqno({\rm MT}_p)
 $$
 \end{thm}

$({\rm MT}_p)$ is an immediate consequence of the noncommutative
Burkholder-Gundy inequalities  below. It was first proved in
\cite{px-BG} for the tracial case, and then extended to the
general case in \cite{jx-burk}. $({\rm MT}_p)$ fails, of course,
for $p=1$. Randrianantoanina \cite{ran-mtrans} proved, however,
the weak type (1,1) substitute for $p=1$ in the tracial case,
which we recall as follows. Let $\f$ be tracial and assume that
$L_p(\M)$ is the usual tracial $L_p$-space associated with $(\M,
\f)$. Let $x$ be a finite $L_1$-martingale. Then
 $$
 \f\big(\un_{(\l,\;\8)}(|\sum_{n\ge1} \e_{n} dx_{n}|)\big)
 \le \kappa_{1}\, \frac{\| x \|_1}\l,\quad \forall\; \l>0, \;\forall\;\e_{n}
 = \pm 1. \leqno({\rm MT}_1)
 $$
This weak type $(1,1)$ inequality gives an alternate simple proof
of $({\rm MT}_p)$ in the tracial case. Indeed, note that $({\rm
MT}_p)$ is trivial for $p=2$ (with $\kappa_2=1$) by virtue of the
orthogonality of martingale differences in $L_2(\M)$. Then
interpolating this trivial case with  $({\rm MT}_1)$ via the
Marcinkiewicz interpolation theorem, we get $({\rm MT}_p)$ for
$1<p<2$. Finally, duality yields the case $2<p<\8$.

We also emphasize that Randrianantoanina's theorem provides a key
to the problem of finding the optimal orders of the best constants
in various noncommutative martingale inequalities. See Corollary
\ref{const mart} below and the discussion following it.

\medskip\n{\bf Burkholder-Gundy inequalities.}  To state these
inequalities, we need to recall the Hardy spaces of noncommutative
martingales introduced in \cite{px-BG}. Let $1 \leq p< \infty$. We
define $\H^c_p(\M)$ to be the space of all martingales $x =
(x_{n})_{n}$  in $L_p(\M)$ such that $dx\in L_p(\M;\el_2^c)$. We
equip $\H^c_p(\M)$ with the norm
 $$\|x\|_{\H_p^c}=\big\|\big(\sum_{n\ge 1}
 |d x_{n}|^{2}\big)^{1/2}\big\|_{p}\,.$$
Similarly,  $\H^r_p(\M)$ is defined to be the space of all
 $L_p$-martingales $x$ such that $x^*\in\H_p^c(\M)$, equipped with
 the norm $\|x\|_{\H^r_p}=\|x^*\|_{\H_p^c}$.
 Finally, set
 $$\H_p(\M)=\H^c_p(\M)\cap \H^r_p(\M)\;\mbox{for}\; p\ge 2
 \quad\mbox{and}\quad
 \H_p(\M)=\H_p^c(\M)+ \H_p^r(\M)\;\mbox{for}\; p<2\,,$$
equipped with the intersection and sum norms, respectively.

\begin{thm}
Let $1 < p < \infty.$ Then for all finite noncommutative
$L_p$-martingales $x = (x_{n})$,
 $$\a^{-1}_p\,\|x\|_{\H_p} \le \|x\|_{p}\le \b_{p} \,
 \|x\|_{\H_p}. \leqno({\rm BG}_p)$$
 \end{thm}

These inequalities were  proved in \cite{px-BG} for the tracial
case and in \cite{jx-burk} for the general case.  The second
inequality in $({\rm BG}_p)$ remains valid for $p=1$, while the
first one has a weak type $(1,1)$ substitute for $p=1$ in the
tracial case (see \cite{ran-weak}). Since the norm
$\|\cdot\|_{\H_p}$ is unconditional on martingale difference
sequences, $({\rm BG}_p)$ immediately implies $({\rm MT}_p)$.
Conversely, by virtue of the noncommutative Khintchine
inequalities  $({\rm K}_p)$, $({\rm MT}_p)$ implies $({\rm BG}_p)$
in the case $p\ge 2$. For $p<2$ we further need the noncommutative
Stein inequality, which is the following statement.

\begin{thm}
  Let $1<p<\8$. Then for all finite sequences
$(a_{n})_{n}$ in $L_p(\M)$,
 $$\big\|\big(\sum_n|\E_n(a_n)|^{2}\big)^{1/2}\big\|_p
 \leq\gamma_{p}\,\big\|\big(\sum_n|a_n|^{2}\big)^{1/2}\big\|_p.
 \leqno({\rm S}_p)$$
 \end{thm}

This result was proved in \cite{px-BG} for the tracial case and in
\cite{jx-burk} for the general noncommutative $L_p$-spaces. We
emphasize that $({\rm S}_p)$ often plays an important role when
dealing with noncommutative martingales. In the tracial case,
$({\rm MT}_1)$ implies a weak type $(1,1)$ substitute of $({\rm
S}_p)$ for $p=1$, which, together with interpolation, provides
another proof of $({\rm S}_p)$.

\medskip\n {\bf Burkholder inequalities.} These inequalities are
closely related with  $({\rm BG}_p)$. To state them we need  more
notation. Let $1\le p<\8$. Let $x=(x_n)_{n\ge1}$ be a finite
martingale in $\M_aD^{1/p}$ (recalling that $\M_a$ denotes the
family of analytic elements of $\M$). Define (with $\E_0=\E_1$)
 $$\|x\|_{h_p^c}=\big\|\big(\sum_{n\ge1}
 \E_{n-1}(|dx_{n}|^2)\big)^{1/2}\big\|_p\,.$$
This defines a norm on the vector space of all finite martingales
in $\M_aD^{1/p}$. The corresponding completion is denoted by
$h_p^c(\M)$. Similarly, we define $h_p^r(\M)$ by passing to
adjoints. Finally, $h_p^d(\M)$ denotes the subspace of
$\el_p(L_p(\M))$ consisting of martingale differences. Then we
define the conditioned version $h_p(\M)$ of $\H_p(\M)$: for $2\leq
p<\8$,
 $$h_p(\M) =h_p^d(\M)\cap h_p^c(\M)\cap h_p^r(\M)\,$$
and for $1\leq p<2$,
 $$h_p(\M)=h_p^d(\M)+h_p^c(\M)+ h_p^r(\M)\,.$$
$h_p(\M)$ is equipped with the intersection or sum norm according
to $p\ge 2$ or $p<2$. We refer to \cite{jx-burk} for more
information on these spaces. The following theorem gives the
noncommutative Burkholder inequalities from \cite{jx-burk}.

\begin{thm}
 Let $1<p<\infty$. Then an $L_p$-martingale $x$
is bounded in $L_p(\M)$ iff $x$ belongs to $h_p(\M)$; moreover, if
this is the case, then
 $$\eta_{p}^{-1}\|x\|_{h_p}\le
 \|x\|_{p}\le\zeta_{p}\,\|x\|_{h_{p}}. \leqno({\rm B}_p)$$
 \end{thm}

Inequalities $({\rm BG}_p)$ and $({\rm B}_p)$ are linked together
through the dual form of the noncommutative Doob maximal
inequality (see Remark~\ref{dual doob} below). The second
inequality of $({\rm B}_p)$ remains true for $p=1$. In the tracial
case, Randrianantoanina \cite{ran-cond} obtained a weak type
$(1,1)$ substitute for the first inequality with $p=1$.

\medskip\n {\bf Rosenthal inequalities.} We first recall the
classical Rosenthal inequalities. Let $(\O, P)$ be a probability
space and let $(f_n)$ be an independent sequence of random
variables in $L_p(\O)$ with $2\le p<\8$. Then
 $$\big\|\sum_n f_n\big\|_p\sim \big(\sum_n\|f_n\|_p^p\big)^{1/p}
 + \big(\sum_n\|f_n\|_2^2\big)^{1/2}\,.$$
To state the noncommutative analogue of this we need to define
independence in the noncommutative setting. In contrast with the
classical case, there now exist several different notions of
independence. The following definition is general enough to
embrace most existing independences.

Let $\N$ and $\A_n$ be $\s^\f$-invariant von Neumann subalgebras
of $\M$ such that $\N\subset\A_n$ for every $n$. The sequence
$(\A_n)$ can be finite. Note that the $\s^\f$-invariance of $\N$
implies that there exists a normal faithful conditional
expectation $\E_\N\,:\,\M\to\N$ preserving the state $\f$.  We say
that $(\A_n)$ are independent over $\N$ or with respect to $\E_\N$
if for every $n$, $\E_\N(ab)=\E_\N(a)\E_\N(b)$ holds for all $a\in
\A_n$ and $b$ in the von Neumann subalgebra generated by
$(\A_k)_{k\neq n}$. A sequence $(a_n)\subset L_p(\M)$ is said to
be independent with respect to $\E_\N$ if there exist $\A_n$ such
that $a_n\in L_p(\A_n)$ and $(\A_n)$ is  independent with respect
to $\E_\N$.

Note that if $(\A_n)$ are independent over $\N$ and $a_n\in
L_p(\A_n)$ with $\E_\N(a_n)=0$, then $(a_n)$ is a martingale
difference sequence relative to the filtration $\big(VN(\A_1, ...,
\A_{n})\big)_{n\ge1}$, where $VN(\A_1, ..., \A_{n})$ is the von
Neumann subalgebra generated by $\A_1, ..., \A_{n}$. We refer to
\cite{jx-ros} for more information.

\begin{thm}
 Let $2\le p<\infty$ and $(a_n)\in L_p(\M)$ be a finite
independent sequence such that $\E_\N(a_n)=0$. Then
 $$
 \frac12\,\big\|(a_n)\big\|_{\R_p}
 \le \big\|\sum_{n} a_n \big\|_p
  \le\nu_p\,\big\|(a_n)\big\|_{\R_p} \,,\leqno({\rm R}_p)$$
where
 $$\big\|(a_n)\big\|_{\R_p}
 =\max\big\{\big(\sum_n\|a_n\|_p^p\big)^{1/p}\,,\;
 \big\|\big(\sum_n\E_\N(|a_n|^2)\big)^{1/2}\big\|_p\,,\;
 \big\|\big(\sum_n\E_\N(|a_n^*|^2)\big)^{1/2}\big\|_p\big\}.$$
 \end{thm}

This result is proved in \cite{jx-ros}. Dualizing the inequality
above, we get a similar one for $1<p<2$ (see \cite{jx-ros} for
more related results).


\subsection{Reduction}


In this subsection we show that all preceding inequalities can be
reduced to the tracial case.

 \begin{thm}\label{red mart}
  If the preceding
inequalities $({\rm MT}_p)$, $({\rm BG}_p)$, $({\rm S}_p)$, $({\rm
B}_p)$ and $({\rm R}_p)$ hold in the tracial case, then they also
hold in the general case with the same relevant best constants.
 \end{thm}

As already noted before, the Khintchine inequalities $({\rm K}_p)$
can also be  reduced to the tracial case. The situation for $({\rm
K}_p)$ is much simpler than all other martingale inequalities.
This is why $({\rm K}_p)$ is not included in the preceding
statement.

As a corollary, we deduce that the estimates on the best constants
in these inequalities for the tracial case are also valid for the
general case. The notation $A\approx B$ means that $C_1 A\le B\le
C_2A$ for two universal positive constants $C_1$ and $C_2$.

 \begin{cor}\label{const mart}
  The  best constants in
$({\rm MT}_p)$, $({\rm BG}_p)$, $({\rm S}_p)$, $({\rm B}_p)$ and
$({\rm R}_p)$ are estimated as follows:
 \begin{enumerate}[\rm(i)]
 \item  $\kappa_p\approx p$ as $p\to\8$;
 \item $\a_p\approx p$ as $p \to \infty$ and
 $\a_p\approx (p-1)^{-1}$ as $p \to 1$;
 \item $\b_p\approx p$ as $p\to\8$ and
 $\b_p\approx 1$ as $p \to 1$;
 \item $\gamma_p\approx p$ as $p \to \infty$;
 \item  $\eta_p\approx (p-1)^{-1}$ as $p \to 1$ and
 $\eta_p\approx p$ as $p \to \infty$;
 \item  $\zeta_p\approx 1$ as $p \to 1$ and
 $\zeta_p\le C p$ for $2\le p<\8$;
 \item $\nu_p\le Cp$ for $2\le p<\8$.
 \end{enumerate}
 \end{cor}

Using his weak type (1,1) inequality for martingale transforms,
Randrianantoanina \cite{ran-mtrans} proved $\kappa_p\le C\,p$ as
$p\to\8$. This is the optimal order of $\kappa_p$ for it is
already so in the commutative case. The optimal order of $\a_p$ as
$p\to\8$ was established in \cite{jx-const} and is the square of
its commutative counterpart. On the other hand, left open in
\cite{jx-const}, the case of $p$ close to 1 was solved  in
\cite{ran-weak}. This order is the same as that in the commutative
case.  The optimal orders of $\b_p$ were determined in
\cite{jx-burk} for $p\le 2$ and in \cite{jx-const} for $p\ge2$.
They are the same as their commutative counterparts. Note that
Pisier \cite{pis-porth} also showed that $\beta_{p}={\rm O}(p)$
for even integers $p$. The estimate $\gamma_p\le C\,p$ was
obtained in \cite{ran-mtrans}. It was proved in \cite{jx-const}
that this is optimal. This optimal order of $\gamma_p$ is the
square of that in the commutative case. The estimates on $\eta_p$
and $\zeta_p$ mainly come from \cite{ran-cond}, although some
partial results already appeared  in \cite{jx-burk} and
\cite{jx-const}. The estimate on $\nu_p$ was obtained in
\cite{jx-ros}. Thus, at the time of this writing, the only
undetermined optimal orders are on $\zeta_p$ and $\nu_p$ as
$p\to\8$. Recall that in the commutative case, both orders are
${\rm O}(p/\log p)$.

\smallskip

The rest of this section is devoted to the proof of Theorem
\ref{red mart}. The idea is simply as follows. We first lift
martingales from $L_p(\M)$ to $L_p(\R(\M))$ and then pull them
back. Remember that we embed $\M$ into the crossed product
$\R(\M)$, and the filtration $(\M_n)_n$ of $\M$ into the
filtration $(\R(\M_n))_n$ of $\R(\M)$. Thus any martingale $x$
relative to $(\M_n)_n$ is also viewed as a martingale relative to
$(\R(\M_n))_n$. This is the lifting procedure. Once we are in
$\R(\M)$, in order to apply  results in the tracial case, we use
the conditional expectations $\Phi_m$ to go to the finite algebras
$\R_m(\M)$. Finally, we return back to $\M$ via the conditional
expectation $\Phi$. This roundabout procedure is almost
transparent for $({\rm MT}_p)$ and $({\rm S}_p)$.  For the others
a little more effort is needed.

\medskip\n {\it Proof of Theorem \ref{red mart}.} This proof is
divided into several parts according to the inequalities in
consideration.

\smallskip\n (i) {\it Reduction of $({\rm MT}_p)$.}
Fix a finite martingale $x$ in $L_p(\M)$ relative to $(\M_n)_n$.
Lifting $x$ to $L_p(\R(\M))$, we consider $x$ as a martingale
relative to $(\R(\M_n))_n$. Then using the conditional expectation
$\Phi_m$ we compress $x$ into a martingale in $L_p(\R_m(\M))$ for
every fixed $m$. More precisely, put
 $$x^{(m)}=\Phi_m(x)=(\Phi_m(x_n))_{n\ge 1}\,,\quad m\in\nat.$$
Then $x^{(m)}$ is a martingale in $L_p(\R_m(\M))$ relative to
$(\R_m(\M_n))_n$. Now $\R_m(\M)$ admits a normal faithful finite
trace $\f_m$. By the construction of $\f_m$ in
section~\ref{Reduction} and \eqref{f-preserving cond}, $\f_m$ is
invariant under $\wh\E_n$ for all $n\ge 1$. On the other hand, the
martingale structure determined by $(\R_m(\M_n))_n$ in
$L_p(\R_m(\M))$ coincides with that in the tracial noncommutative
$L_p$-space of $\R_m(\M)$. Thus using $({\rm MT}_p)$ in the
tracial case, we have
 $$\big\|\sum_{n\ge1}\e_n d x^{(m)}_n\big\|_p
 \leq\kappa_{p}\, \|x\|_{p}\,,\quad \forall\,\e_{n}=\pm 1.$$
Let $y=\sum_{n\ge1} \e_{n}dx_{n}$. By \eqref{cond hat cond hat},
the sum on the left-hand side is equal to $\Phi_m(y)$.  Now
$(\Phi_m(y))_m$ is a martingale in $L_p(\R(\M))$ with respect to
$(\R_m(\M))_m$. Hence, $\Phi_m(y)\to y$ in $L_p(\R(\M))$ as
$m\to\8$ (see Remark \ref{mart conv}). Therefore,
  $$\big\|\sum_{n\ge1}\e_n d x_{n}\big\|_p
 \leq\kappa_{p}\, \|x\|_p\,.$$
This is $({\rm MT}_p)$ in the general case.

\smallskip\n (ii) {\it Reduction of $({\rm S}_p)$.} This is also easy.
The main point is the following.

\begin{lem}\label{column conv}
 Let $1\le p<\8$. Let $a=(a_n)\in L_p(\M; \el_2^c)$.
Considering $a$ as an element in $L_p(\R(\M);\el^c_2)$, we set
 $$ a^{(m)}=\Phi_m(a)=(\Phi_m(a_n))_{n\ge 1}\,.$$
Then $a^{(m)}\to a$ in $L_p(\R(\M);\el_2^c)$ as $m\to\8$. A
similar statement holds for the row space.
 \end{lem}

\pf We regard $L_p(\M; \el_2^c)$ as  the column subspace of
$L_p(\M\bar\ot B(\el_2))$. Then $L_p(\R(\M); \el_2^c)$ is the
column subspace of $L_p(\R(\M)\bar\ot B(\el_2))$. Note that
 $$a^{(m)}=\Phi_m\ot{\rm id}_{B(\el_2)}(a).$$
Since $(\Phi_m)_m$ is an increasing sequence of conditional
expectations, the desired result immediately follows from the
martingale mean convergence  in Remark \ref{mart conv}.\cqd

\smallskip

Now as for $({\rm MT}_p)$, it is easy to see why one needs only to
consider $({\rm S}_p)$ in the tracial case. Indeed, fix a finite
sequence $a=(a_n)\subset L_p(\M)\subset L_p(\R(\M))$. Then
$\Phi_m(a)\subset L_p(\R_m(\M))$. Applying $({\rm S}_p)$ in the
tracial case, we get
 $$\big\|\big(\E_n(\Phi_m(a_n))\big)_n\big\|_{L_p(\R_m(\M);\el_2^c)}
 \le\gamma_p\,\|\Phi_m(a)\|_{L_p(\R_m(\M);\el_2^c)}\,,\quad m\in\nat.$$
By \eqref{cond hat cond hat},
$\E_n(\Phi_m(a_n))=\Phi_m(\E_n(a_n))$ for all $m,\;
 n\in\nat$. Therefore,
 $$\big\|\big(\Phi_m(\E_n(a_n))\big)_n\big\|_{L_p(\R(\M);\el_2^c)}
 \le\gamma_p\|\Phi_m(a)\|_{L_p(\R(\M);\el_2^c)}\,.$$
It remains to apply Lemma \ref{column conv} to conclude the
reduction argument on $({\rm S}_p)$.

\medskip\n (iii) {\it Reduction of $({\rm BG}_p)$.} In the
following the Hardy spaces on $\R(\M)$ and $\R_m(\M)$ are relative
to the filtrations  $(\R(\M_n))_n$ and $(\R_m(\M_n))_n$,
respectively.

\begin{lem}\label{hardy compl}
 Let $1\le p<\8$. Then $\H_p(\M)$ and $\H_p(\R_m(\M))$ are
$1$-complemented isometric subspaces of $\H_p(\R(\M))$.
 \end{lem}

\pf We consider only the part on $\H_p(\M)$, that on
$\H_p(\R_m(\M))$ being dealt with by the same arguments. It is
trivial that $\H_p^c(\M)$ (resp. $\H_p^r(\M)$) is an isometric
subspace of $\H_p^c(\R(\M))$ (resp. $\H_p^r(\R(\M))$). Let $\wh
x=(\wh x_n)_n$ be a martingale relative to $(\R(\M_n))_n$. Put
$\Phi(\wh x)= (\Phi(\wh x_n))_n$. By \eqref{cond hat cond hat},
$\Phi(\wh x)$ is a martingale relative to $(\M_n)_n$. Note that
the difference sequence of $\Phi(\wh x)$ is given by
 $(\Phi(d\wh x_n))_n$.
Then using the tensor product argument in the proof of Lemma
\ref{column conv},  one sees that the map $\wh x\mapsto\Phi(\wh
x)$ defines a contractive projection from $\H_p^c(\R(\M))$ (resp.
$\H_p^r(\R(\M))$) onto $\H_p^c(\M)$ (resp. $\H_p^r(\M)$).
Therefore, we immediately deduce the assertion on $\H_p(\M)$ in
the case $p\ge2$. For the case $p<2$ we need only to check that
$\H_p(\M)$ is an isometric subspace of $\H_p(\R(\M))$. But this is
an easy consequence of the above projection.\cqd

\begin{lem}\label{hardy conv}
 Let $1\le p<\8$. Let $x\in \H_p(\M)$. Then
$\Phi_m(x)\in\H_p(\R(\M))$ and $\Phi_m(x)\to x$ in $\H_p(\R(\M))$
as $m\to\8$.
 \end{lem}

\pf By Lemma \ref{column conv} we see that the statement above
holds with $\H_p^c(\R(\M))$ or $\H_p^r(\M)$ in place of
$\H_p(\R(\M))$. We then deduce the assertion in the case $p\ge 2$.
The case $p<2$ is proved  by virtue of the density of finite
martingales in $\H_p(\R(\M))$ and the contractivity of $\Phi_m$ on
$\H_p(\R(\M))$. \cqd

\smallskip

Using the previous two lemmas, we easily reduce the general case
of $({\rm BG}_p)$ to the tracial case, as before for $({\rm
S}_p)$. We leave the details to the reader.

\medskip\n (iv) {\it Reduction of $({\rm B}_p)$.} The proof of
this is very similar to the previous one. It suffices to apply the
following lemma, which is the analogue of Lemmas \ref{hardy compl}
and \ref{hardy conv} for the conditioned Hardy spaces.

\begin{lem}
 Let $1\le p<\8$.
 \begin{enumerate}[\rm(i)]
 \item Then $h_p(\M)$ and $h_p(\R_m(\M))$ are
$1$-complemented isometric subspaces of $h_p(\R(\M))$.
 \item $\Phi_m$ converges to the identity of $h_p(\R(\M))$
in the point-norm topology.
 \end{enumerate}
 \end{lem}

\pf  We first show that $\Phi$ projects $h_p^k(\R(\M))$
contractively onto $h_p^k(\R(\M))$ for $k\in\{d, c, r\}$. This is
trivial for $h_p^d(\R(\M))$. To deal with $h_p^c(\R(\M))$ we first
recall that $h_p^c(\R(\M))$ can be isometrically viewed as the
column subspace of $L_p(\R(\M)\bar\ot B(\el_2(\nat^2)))$. The map
realizing this is constructed by means of Kasparov's stabilization
theorem for Hilbert C*-modules (see \cite{ju-doob} for more
details). Using this and as in the proof of Lemma \ref{hardy
compl}, we show that $\Phi$ is a contractive projection on
$h_p^c(\R(\M))$. Passing to adjoints, we get the assertion on
$h_p^r(\R(\M))$. We then deduce (i). (ii) is proved similarly as
Lemma \ref{hardy conv}, so we omit the details. \cqd

\medskip\n (v) {\it Reduction of $({\rm R}_p)$.} Let $(\A_n)$ be
a sequence of von Neumann subalgebras of $\M$ independent over
$\N$. We use the notation introduced in subsection \ref{The
framework}. In particular, $\R(\A_n)$ and $\R(\N)$ are
$\s^{\wh\f}$-invariant von Neumann subalgebras of $\R(\M)$.

\begin{lem}
 The algebras $\R(\A_n)$ are independent over $\R(\N)$.
 \end{lem}

\pf We must show that
 \beq\label{l-indep}
 \E_{\R(\N)}(\hat a\hat b)=\E_{\R(\N)}(\hat a)\,\E_{\R(\N)}(\hat b)
 \eeq
for $\hat a\in \R(\A_n)$ and $\hat b$ in the subalgebra generated
by the $\R(\A_k)$, $k\neq n$. Let $\B$ be the von Neumann
subalgebra of $\M$ generated by $(\A_k)_{k\neq n}$. Note that $\B$
is $\s$-invariant. Using \eqref{com pi-l}, we see that the
subalgebra generated by $(\R(\A_k))_{k\neq n}$ is equal to the
crossed product $\R(\B)=\B\rtimes_{\s} G$.  On the other hand,
$\E_{\R(\N)}=\wh{\E_\N}$ by Theorem \ref{extension crossed}. Now
let $a\in\A_n$, $b\in\B$ and $g, h\in G$. Then by \eqref{com
pi-l}, \eqref{def of T hat} and the independence of $(\A_n)_n$ we
find
 \be
 \E_{\R(\N)}\big(\l(g)a\,\l(h)b\big)
 &=&\wh{\E_\N}\big(\l(g)\l(h)\,\s_{-h}(a)b\big)\\
 &=&\l(g)\l(h)\,\E_\N\big(\s_{-h}(a)b\big)\\
 &=&\l(g)\l(h)\,\E_\N\big(\s_{-h}(a)\big)\,\E_\N(b)\\
 &=&\l(g)\l(h)\,\s_{-h}\big(\E_\N(a)\big)\,\E_\N(b)\\
 &=&\l(g)\E_\N(a)\,\l(h)\E_\N(b)\\
 &=&\E_{\R(\N)}\big(\l(g)a\big)\, \E_{\R(\N)}\big(\l(h)b\big).
 \ee
This shows \eqref{l-indep} for $\hat a=\l(g)a$ and $\hat
b=\l(h)b$. Since $\R(\A_n)$ (resp. $\R(\B)$) is the w*-closure of
all finite sums on $\l(g)a$ with $a\in\A_n, g\in G$ (resp.
$\l(h)b$ with $b\in\B, h\in G$), the normality of $\E_{\R(\N)}$
implies \eqref{l-indep} for general $\hat a$ and $\hat b$. \cqd

\smallskip

It is then an easy exercice to reduce the general case of $({\rm
R}_p)$ to the tracial one. We omit the details. Therefore, the
proof of Theorem \ref{red mart} is complete. \cqd


\section{Applications to noncommutative maximal inequalities}
 \label{Applications to noncommutative maximal inequalities}


We pursue our applications of Theorem \ref{red s-finite} to
noncommutative inequalities. We now consider the noncommutative
maximal martingale and ergodic inequalities. We should emphasize
that because of the failure of the noncommutative analogue of the
usual pointwise maximal function of a sequence of functions, these
maximal inequalities are much subtler than those in the
commutative case. This failure also explains why we are forced to
work systematically with the vector-valued spaces
$L_p(\M;\el_\8)$. Namely, instead of a noncommutative pointwise
maximal function (which does not exist now), we work with the
noncommutative analogue of the usual maximal $L_p$-norm. In this
section $\M$ again denotes a von Neumann algebra equipped with a
normal faithful state $\f$.


\subsection{Positive maps on $L_p(\M;\el_\8)$}


We first recall the definition of the spaces $L_p(\M;\ell_\infty)$
and $L_p(\M;\ell_1)$ from \cite{ju-doob}. Let $1\le p\le\8$. A
sequence $(x_n)$ in $L_p(\M)$ belongs to $L_p(\M;\ell_\8)$ iff
$(x_n)$ admits a factorization $x_n=ay_nb$ with $a, b\in
L_{2p}(\M)$ and $(y_n)\in\ell_\8(L_\8(\M))$. The norm of $(x_n)$
is then defined as
 $$
  \|(x_n)\|_{L_p(\M;\ell_\8)}
 =\inf_{x_n=ay_nb}\, \|a\|_{2p}\,
 \sup_n\|y_n\|_\8\, \|b\|_{2p}\,.
 $$
On the other hand,  $L_p(\M;\ell_1)$ is defined to be the space of
all sequences $(x_n)\subset L_p(\M)$ for which there exist
$a_{nk}, b_{nk}\in L_{2p}(\M)$ such that
 \[x_n=\sum_ka_{nk}^*b_{nk}\,, \quad \forall\;n\ge1.\]
$L_p(\M;\ell_1)$ is equipped with the norm
 \[ \|(x_n)\|_{L_p(\M;\ell_1)}
 =\inf_{x_n=\sum_{k} a_{nk}^*b_{nk}}\,
 \big\|\sum_{n,k} a_{nk}^*a_{nk}\big\|_p^{1/2}\,
 \big\|\sum_{n,k} b_{nk}^*b_{nk}\big\|_p^{1/2} \,.\]
We refer to \cite{ju-doob} for more information (see also
\cite{jx-erg}). Let us note that if $\M$ is injective, the
definition above is a special case of Pisier's vector-valued
noncommutative $L_p$-space theory \cite{pis-ast}. The following
remark is easy to check.

\begin{rk}\label{positive Lplinfty}
 Let $(x_n)\subset L_p^+(\M)$.
 \begin{enumerate}[(i)]
 \item $(x_n)\in L_p(\M;\el_\8)$ iff there exists $a\in
L_p^+(\M)$ such that $x_n\le a$ for all $n\in\nat$. In this case,
we have
 $$\|(x_n)\|_{L_p(\M;\el_\8)}=\inf\big\{\|a\|_p\;:\; a\in
 L_p^+(\M)\;\mbox{s.t.}\; x_n\le a,\;\forall\;n\big\}.$$
 \item $(x_n)\in L_p(\M;\el_1)$ iff $\sum_nx_n\in L_p(\M)$.
If this is the case, then
 $$\|(x_n)\|_{L_p(\M;\el_1)}=\big\|\sum_{n\ge1}x_n\big\|_p\,.$$
 \end{enumerate}
  \end{rk}

We adopt the convention of \cite{jx-erg} that the norm
$\|(x_n)\|_{L_p(\M;\el_\8)}$ is denoted by
$\big\|\sup_n^+x_n\big\|_p\,$. However, we should warn the reader
that $\big\|\sup^+_nx_n\big\|_p$ is just a notation, for
$\sup_nx_n$ does not make any sense in the noncommutative setting.
We should also point out that $L_p(\M;\el_\8)$ and $L_p(\M;\el_1)$
are not closed under absolute value. In particular,
$\big\|\sup^+_nx_n\big\|_p\neq \big\|\sup^+_n|x_n|\big\|_p$ in
general.

\smallskip

One fundamental result on these vector-valued spaces is the
duality between $L_p(\M;\el_1)$ and $L_{p'}(\M;\el_\8)$ for $1\le
p<\8$ established in \cite{ju-doob}, where $p'$ denotes the
conjugate index of $p$. Namely, we have
 \beq\label{Lplinfty-dual1}
 L_p(\M;\el_1)^*=L_{p'}(\M;\el_\8)\quad\mbox{isometrically}
 \eeq
via the duality bracket
 $$\la (x_n),\;(y_n)\ra=\sum_{n\ge1}\tr(x_ny_n), \quad
 (x_n)\in L_p(\M;\el_1),\; (y_n)\in L_{p'}(\M;\el_\8).$$
Moreover, if $L_p(\M;\el_q^n)$ denotes the $\el_q^n$-valued
analogues of these spaces ($q=1,\8$), then
 \beq\label{Lplinfty-dual2}
 L_p(\M;\el_\8^n)^*=L_{p'}(\M;\el_1^n)\quad\mbox{isometrically}
 \eeq
(see \cite{jx-erg}). Using these duality identities and the
following formula,
 \beq\label{Lplinfty1}
 \big\|\,{\sup_n}^+x_n\big\|_p=\sup_{n\ge 1}
 \big\|\mathop{{\sup}^+}_{1\le k\le n}x_k\big\|_p
 \eeq
from \cite{jx-erg}, we deduce the following.

\begin{rk}\label{sup via sum}
 Let $(x_n)\in L_p(\M;\el_\8)$, $1\le p\le\8$. Then
 $$\big\|\,{\sup_n}^+x_n\big\|_p=\sup\big\{
 \big|\sum_{n\ge1}\tr(x_ny_n)\big|\;:\; (y_n)\in
 L_{p'}(\M;\el_1),\; \big\|(y_n)\big\|_{L_{p'}(\M;\el_1)}\le
 1\big\}.$$
Moreover, if $(x_n)$ is positive, the supremum above can be
restricted to positive $(y_n)$ too.
  \end{rk}

Recall that a map $T\,:\, L_p(\M)\to L_p(\M)$ is $n$-positive if
$T_n\,:\, L_p(\ma_n\ot\M)\to L_p(\ma_n\ot\M)$ is positive, where
$\ma_n$ denotes the full algebra of $n\times n$ matrices and where
 $$T_n\big((x_{ij})_{1\le i,j\le n}\big)
 =\big(T(x_{ij})\big)_{1\le i,j\le n}\,. $$
Here we have viewed, as usual, the elements of $L_p(\ma_n\ot\M)$
as $n\times n$-matrices with entries in $L_p(\M)$. $T$ is
completely positive if $T$ is $n$-positive for every $n$.

\begin{prop}\label{positive map on Lpelinfty}
 Let $1\le p<\8$ and $T\,:\,L_p(\M)\to L_p(\M)$ be
a positive bounded map. Then $T\ot{\rm id}_{\el_\8}$ is bounded on
$L_p(\M;\el_\8)$ and of norm $\le 8\|T\|$. If in addition $T$
 is $2$-positive, then $\|T\ot{\rm id}_{\el_\8}\|=\|T\|$.
  \end{prop}

\pf By \eqref{Lplinfty-dual1}, \eqref{Lplinfty-dual2} and
\eqref{Lplinfty1}, it suffices to show that $T^*$ extends to a
bounded map on $L_{p'}(\M;\el_1)$. However, by \cite{jx-erg},
every element in the unit ball of $L_{p'}(\M;\el_1)$ is a sum of
$8$ elements in the same ball. Thus we need only to consider a
positive $x=(x_n)\in L_{p'}(\M;\el_1)$. Then $T^*(x_n)\ge 0$ too.
Therefore, by Remark \ref{positive Lplinfty},
 $$\big\|T^*\ot{\rm id}_{\el_1}(x)\big\|_{L_{p'}(\M;\el_1)}=
 \big\|\sum_nT^*(x_n)\big\|_{p'}\le \|T\|\,
 \big\|\sum_nx_n\big\|_{p'}\,.$$
It thus follows that $\|T^*\ot{\rm id}_{\el_1}\|\le 8\|T\|$. This
yields the first assertion.

The proof of the second assertion needs more effort. Recall that
$D$ denotes the density operator in $L_1(\M)$ of the distinguished
state $\f$. Let
 $$\wt\M=\ma_2\ot\M\quad\mbox{and}\quad \wt\f=\Tr\ot\f\,,$$
where $\Tr$ is the usual trace on $\ma_2$. Then the density of
$\wt\f$ is the matrix
 $$\wt D=\begin{pmatrix}
 D&0\\
 0&D
 \end{pmatrix}\,.$$
Let $a, b\in L_{2p}(\M)$ be positive operators and $\e>0$. Put
 $$d=\begin{pmatrix}
 a&0\\
 0&b
 \end{pmatrix}\quad\mbox{and}\quad d_\e=T_2(d^2) +\e\,{\wt D}^{1/p}\,.$$
Note that $d_\e\in L_p(\wt\M)$ and $d_\e$ is an injective positive
operator affiliated with the crossed product
$\wt\M\rtimes_{\s^{\wt\f}}\real$. Define
 $$\wt T(\tilde z)=d_\e^{-1/2}\,T_2(d\tilde zd)\,
 d_\e^{-1/2}\,,
 \quad \tilde z\in\wt\M\,.$$
The $2$-positivity of $T$ implies that
 $$0\le T_2(d\tilde zd)\le\|\tilde z\|_\8\, T_2(d^2)
 \le\|\tilde z\|_\8\, d_\e,\quad \forall\;
 \tilde z\in\wt\M_+\,.$$
It follows that $\wt T(\tilde z)$ is bounded for every $\tilde
z\ge0$; so $\wt T(\tilde z)\in \wt\M\rtimes_{\s^{\wt\f}}\real$. On
the other hand, it is easy to see that $\wt T(\tilde z)$ is
invariant under the dual action of $\s^{\wt\f}$. Therefore, $\wt
T(\tilde z)\in \wt\M$ for $\tilde z\in\wt\M_+$. Decomposing every
element into a linear combination of 4 positive elements, we then
deduce that $\wt T$ is a well-defined map on $\wt\M$. Moreover,
$\wt T$ is positive. Thus
 $$\|\wt T\|=\|\wt T(1)\|\le1.$$
For $z\in \M$ define
 $$\tilde z=\begin{pmatrix}
 0 &z\\ 0&0 \end{pmatrix}\,.$$
Then
 $$\wt T(\tilde z)=\begin{pmatrix}
 0 & a_\e^{-1/2}\, T(azb)\,b_\e^{-1/2}
 \\ 0&0 \end{pmatrix}\,,$$
where $a_\e=T(a^2)+\e\, D^{1/p}$ and $b_\e=T(b^2)+\e\, D^{1/p}$.
Therefore, the map $T'\,:\, \M\to\M$ defined by
$T'(z)=a_\e^{-1/2}\, T(azb)\,b_\e^{-1/2}$ is a contraction.

Now it is easy to finish the proof. Indeed, let $x=(x_n)\in
L_p(\M;\el_\8)$ and take a factorization $x_n=ay_nb$ of $x$ with
$a, b\in L_{2p}(\M)$ and $(y_n)\in \el_\8(L_\8(\M))$. By polar
decomposition we can assume $a, b\ge0$. Using the map $T'$ above
we get a factorization of $(T(x_n))$ as follows:
 $$T(x_n)=a_\e^{1/2}\, T'(y_n)\, b_\e^{1/2}
 \;{\mathop =^{\rm def}}\; a'y'_nb'\,.$$
We have
 $$\|a'\|_{2p}\le\|T\|^{1/2}(\|a\|_{2p} +{\rm O}(\e)),\quad
 \|b'\|_{2p}\le\|T\|^{1/2}(\|b\|_{2p} +{\rm O}(\e)),\quad
 \|y_n'\|_\8\le\|y_n\|_\8\,.$$
It follows that $(T(x_n))\in L_p(\M;\el_\8)$ and
 $$\big\|\,{\sup_n}^+T(x_n)\big\|_p\le \|T\|
 \big(\big\|\,{\sup_n}^+x_n\big\|_p
 +{\rm O}(\e)\big).$$
This implies the desired assertion. \cqd

\smallskip

Let $\N\subset\M$ be a $\s^\f$-invariant von Neumann subalgebra
and let $\E\,:\,\M\to\N$ be the normal conditional expectation
preserving the state $\f$. Under this assumption we know that
$L_p(\N)$ is naturally identified as a subspace of $L_p(\M)$ and
$\E$ extends to a contractive projection from $L_p(\M)$ onto
$L_p(\N)$ (see \cite[Lemma~1.2]{jx-burk} and
Example~\ref{extension Lp cond}). Moreover, $\E$ is completely
positive. Applying the previous proposition to this situation we
immediately get  the following.

\begin{cor}\label{condi exp on Lplinfty}
 The inclusion $L_p(\N;\el_\8)\subset L_p(\M; \el_\8)$ is
isometric and $\E\ot{\rm id}_{\el_\8}$ defines a contractive
projection from $L_p(\M; \el_\8)$ onto $L_p(\N; \el_\8)$.
 \end{cor}


\subsection{Doob maximal inequality}


In this subsection we keep the framework introduced in subsection
\ref{The framework}. In particular, $(\M_n)_n$ is a filtration of
von Neumann subalgebras of $\M$ with associated conditional
expectations $(\E_n)_n$. The following is the noncommutative Doob
maximal inequality.

\begin{thm}
Let $1<p\le\infty$. Let $x=(x_n)$ be a bounded $L_p$-martingale
with respect to $(\M_n)$. Then
  $$\big\|\,{\sup_n}^+x_n\big\|_p\le \d_p\|x\|_p\,.
  \leqno({\rm D}_p)$$
Moreover, $\d_p\le C\,(p-1)^{-2}$ for a universal constant $C$ and
this estimate is optimal as $p\to1$.
 \end{thm}

For positive martingales, $({\rm D}_p)$ takes the following much
simpler form, from which the reader can recognize the classical
Doob maximal inequality.

 \begin{rk} Let
$x=(x_n)$ be a positive bounded $L_p$-martingale. Then there
exists $a \in L_p^+(\M)$ such that
 $$\|a\|_{p} \leq \delta_{p}^+\, \|x\|_{p}
 \quad\mbox{and }\quad
 x_n\le a,\quad \forall\; n\in\nat. \leqno({\rm D}_{p}^+)$$
 \end{rk}

One can show that $({\rm D}_p)$ and $({\rm D}_p^+)$ are equivalent
and that the two relevant  best constants $\d_p$ and $\d_p^+$ are
equal (see \cite{ju-doob}). However, if we tolerate  a multiple of
constants, we easily see that $({\rm D}_p^+)$ implies $({\rm
D}_p)$ with $\d_p\le 4\d_p^+$.

 \begin{rk}\label{dual doob}
  It is sometimes more convenient to work
with the dual form of $({\rm D}_p)$.  Let $1 \le q < \infty.$ Then
for all finite sequences $(a_{n})_{n}$ of positive elements in
$L_q(\M)$,
 $$\big\|\sum_n\E_n(a_n)\big\|_q\le \d'_q\,\big\|\sum_n
 a_n\big\|_q\,. \leqno({\rm D}'_q)$$
One can show that $({\rm D}_p)$ with $1<p\le\8$ is equivalent to
$({\rm D}'_{p'})$ and $\d_{p} = \delta'_{p'}$. We refer to
\cite{ju-doob} for details.
 \end{rk}

\smallskip

Inequality  $({\rm D}_p)$ was first obtained in \cite{ju-doob},
and the optimal order of $\d_p$ was determined in \cite{jx-const}.
Note that this order is the square of that in the commutative
case. The proof of \cite{ju-doob}  heavily relies upon Hilbert
C*-module theory. An elementary proof for the tracial case was
later given in \cite{jx-erg}. A key ingredient of this second
proof is a noncommutative Marcinkiewicz-type interpolation
theorem. It is also this new proof that gives the optimal order
$\d_p={\rm O}((p-1)^{-2})$ as $p\to1$.

\smallskip

As for the square function type inequalities in the previous
section, it suffices to show $({\rm D}_p)$ in the tracial case.

\begin{prop}
 If $({\rm D}_p)$ holds in the tracial case, so does it in the
general case with the same best constant.
 \end{prop}

\pf Assume that $({\rm D}_p)$ holds in the tracial case with
constant $\d_p$ for $1<p<\8$. Let $x=(x_n)$ be a bounded
$L_p$-martingale relative to $(\M_n)_n$. Set
$x^{(m)}=(\Phi_m(x_n))_n$. Then $x^{(m)}$ is a bounded
$L_p$-martingale relative to $(\R_m(\M_n))_n$. Now $\R_m(\M)$ is a
finite von Neumann algebra. Using $({\rm D}_p)$ in the tracial
case we find
 $$\big\|x^{(m)}\big\|_{L_p(\R_m(\M);\el_\8)}\le\d_p\,
 \big\|x^{(m)}\big\|_p\,.$$
Consequently,
 $$\big\|x^{(m)}\big\|_{L_p(\R(\M);\el_\8)}\le\d_p\,
 \|x\|_p\,.$$
In particular,
 $$\big\|(x_1^{(m)}\,,\cdots, x_n^{(m)})\big\|_{L_p(\R(\M);\el_\8^n)}
 \le\d_p\,\|x\|_p\,,\quad\forall\; m,\; n\in\nat.$$
Note that for a fixed $n$ the norm of $L_p(\R(\M);\el_\8^n)$ is
equivalent to that of  $\el_\8^n(L_p(\R(\M))$. Since $x_k^{(m)}\to
x_k$ in $L_p(\R(\M))$ as $m\to\8$ for every $k\in\nat$, we get
 $$\lim_{m\to\8}
 \big\|(x_1^{(m)}\,,\cdots, x_n^{(m)})\big\|_{L_p(\R(\M);\el_\8^n)}
 =\big\|(x_1\,,\cdots, x_n)\big\|_{L_p(\R(\M);\el_\8^n)}\,.$$
Thus we deduce
 $$\sup_n\big\|(x_1\,,\cdots, x_n)\big\|_{L_p(\R(\M);\el_\8^n)}
 \le\d_p\,\|x\|_p\,.$$
Together with \eqref{Lplinfty1}, this implies that
 $$\big\|x\big\|_{L_p(\R(\M);\el_\8)}
 \le\d_p\,\|x\|_p\,.$$
Recall that $x=(x_n)$ is a martingale in $L_p(\M)\subset
L_p(\R(\M))$. Applying Corollary \ref{condi exp on Lplinfty} to
the conditional expectation $\Phi\,:\, \R(\M)\to\M$, we finally
get
 $$\big\|x\big\|_{L_p(\M;\el_\8)}
 \le\d_p\,\|x\|_p\,.$$
This is the desired inequality $({\rm D}_p)$ in $L_p(\M)$. \cqd


\subsection{Maximal ergodic inequalities}


In this subsection we deal with noncommutative maximal ergodic
inequalities. In the sequel $T$ will be assumed to be a map on
$\M$ satisfying the following properties:
 \begin{itemize}
 \item[(I)] $T$ is a contraction on $\M$;
 \item[(II)] $T$ is completely positive;
 \item[(III)] $\f\circ T\le \f$;
 \item[(IV)] $T\circ\s_t^\f=\s_t^\f\circ T$ for all $t\in\real$.
 \end{itemize}

By Theorem \ref{extension Lp}, $T$ extends to a contraction on
$L_p(\M)$ for any $p\ge1$, the extension  still being denoted by
$T$. Note that by virtue of condition (IV) we can also use the
extension given by Proposition \ref{extension Lp bis}.

\medskip

We now consider the ergodic averages of $T$:
 $$M_n\equiv M_n(T)=\frac{1}{n+1}\,\sum_{k=0}^n T^k\,.$$
The following theorem gives our maximal ergodic inequalities on
$T$.

\begin{thm}\label{haagmax}
 Let $1<p\le\8$. Then
 \beq\label{dsp}
 \big\|\,{\sup_n}^+M_n(x)\big\|_p\le C_p\,\|x\|_p\;,\quad x\in
 L_p(\M).
 \eeq
If in addition $T$  is symmetric in the following sense:
 \beq\label{sym}
 \f(T(y)^*x)=\f(y^*T(x)),\quad\forall\;x,\;y\in\M,
 \eeq
then
 \beq\label{smp}
 \big\|\,{\sup_n}^+ T^n(x)\big\|_p\le C'_p\,\|x\|_p\;,\quad x\in
 L_p(\M).
 \eeq
Here $C_p$ and  $C_p'$ are two positive constants depending only
on $p$.
 \end{thm}

Inequality \eqref{dsp} is the noncommutative analogue of the
classical Dunford-Schwartz maximal ergodic inequality for positive
contractions on $L_p$ (see \cite{dunford} for the classical case),
while \eqref{smp} extends Stein's maximal inequality \cite{st-erg}
(see also \cite{st-lp}) to the noncommutative setting. The
preceding theorem is proved in \cite{jx-erg}. The order of the
constant $C_p$ obtained there is $(p-1)^{-2}$ as $p\to1$. This
order is optimal. However, we do not know the optimal order of
$C'_p$. Recall that in the commutative case both $C_p$ and $C'_p$
are of optimal order $(p-1)^{-1}$ as $p\to1$.

\smallskip

Like the noncommutative Doob maximal inequality, both \eqref{dsp}
and \eqref{smp} can be reduced to the tracial case. Note that if
$\f$ is tracial, condition (IV) is automatically satisfied. The
proof of \eqref{dsp} in the tracial case is based on Yeadon's weak
type $(1, 1)$ ergodic inequality (see \cite{ye1}) and the
noncommutative Marcinkiewicz-type interpolation theorem already
mentioned in the previous subsection. The reduction of \eqref{dsp}
and \eqref{smp} to the tracial case was done in \cite{jx-erg}. For
the convenience of the reader we include its main lines here.

\medskip

\n{\it Reduction of Theorem \ref{haagmax} to the tracial case.} We
will again use Theorem \ref{red s-finite}. Recall that $\R$ is the
crossed product $\M\rtimes_{\s^{\f}}G$ and  $(\R_m)_{m \ge 1}$ is
the filtration  of finite von Neumann subalgebras of $\R$
constructed in Theorem \ref{red s-finite}.

On the other hand, using Theorem \ref{extension crossed} we extend
$T$ to a completely positive map $\wh T$ on $\R$. By that theorem
and Proposition \ref{extension crossed bis}, $\wh T$ satisfies
conditions (I)-(IV) relative to $(\R, \;\wh\f)$. Therefore, $\wh
T$ extends to a contraction on $L_p(\R)$ for every $p\ge1$. Note
that if $T$ satisfies the symmetry condition \eqref{sym}, then
using the arguments in the proof of Theorem \ref{extension
crossed}, we easily check that $\wh T$ is also symmetric relative
to $\wh\f$.

Let $\t_m$ be the trace on $\R_m$ corresponding to the state
$\f_m$ constructed in the proof of Theorem \ref{red s-finite}. By
Proposition \ref{extension crossed bis}, we have $\t_m\circ\wh
T\big|_{\R_m}\le\t_m$. Therefore, $\wh T\big|_{\R_m}$ satisfies
conditions (I)-(IV) relative to $(\R_m,\;\t_m)$ for every $m$.
This will allow us to apply \eqref{dsp} in the tracial case.

Now assume \eqref{dsp} in the tracial case. Fix $1<p<\8$ and $x\in
L_p(\M)$. We consider $x$ as an element in $L_p(\R)$ and then
apply the conditional expectation $\Phi_m$ to it:
 $\displaystyle x_m=\Phi_m(x)\in L_p(\R_m)$.
Applying \eqref{dsp} to $\wh T$ on $\R_m$, we get
 $$\big\|\,{\sup_n}^+ M_n(\wh T)(x_m)\big\|_p\le C_p\,\|x\|_p\;,\quad
 \forall\; m\in\nat.$$
By the martingale convergence theorem in Remark \ref{mart conv},
 $$\lim_{m\to\8}x_m=x\quad \mbox{in}\quad L_p(\R).$$
Consequently,
 $$\lim_{m\to\8}\wh T\,^k(x_m)=x\quad \mbox{in}\quad L_p(\R),
 \quad\forall\; k\ge0.$$
It then follows that
 $$\lim_{m\to\8}\big\|\mathop{{\sup}^+}_{1\le k\le n} M_k(\wh T)(x_m)\big\|_p
 =\big\|\mathop{{\sup}^+}_{1\le k\le n} |M_k(\wh T)(x)|\big\|_p\;.$$
However, since $x\in L_p(\M)$,
 $$M_k(\wh T)(x)=M_k(T)(x).$$
Therefore, we deduce
 $$\big\|\mathop{{\sup}^+}_{1\le k\le n} M_k(T)(x)\big\|_p\le C_p\,\|x\|_p\;,
 \quad\forall\; n\in\nat.$$
Thus, by \eqref{Lplinfty1},
 $$\big\|\,{\sup_n}^+ M_n(T)(x)\big\|_p\le C_p\,\|x\|_p\;.$$
This shows \eqref{dsp} in the general case. The reduction of
\eqref{smp} is dealt with in a similar way.

\begin{rk}
 (i) By discretization we can easily extend Theorem \ref{haagmax}
to the case of semigroups. We refer to \cite{jx-erg} for more
details.

(ii) It is well known that a maximal ergodic theorem usually
yields a corresponding individual ergodic theorem. This is indeed
the case of Theorem \ref{haagmax}. For instance, inequality
\eqref{dsp} implies that $M_n(x)$ converges bilaterally almost
surely in the sense of \cite{ja2} to an element $\wh x$ for every
$x\in L_p(\M)$ with $1<p<\8$. In the case $p=\8$, $M_n(x)$
converges to $\wh x$ almost uniformly in Lance's sense
\cite{lance-erg}. Thus in the latter case, we recover the
individual ergodic theorems of Lance, Conze and
 Dang-Ngoc \cite{cdn} and K\"ummerer \cite{kum}. The interested
 reader is referred to \cite{jx-erg} for more information.
 \end{rk}


\bigskip

 \footnotesize{

\noindent  U.H.: Department of Mathematics and Computer Science,
University
 of Southern Denmark,  Campusvej 55, 5230 Odense M, Denmark.\\
 haagerup@imada.sdu.dk\\

\noindent M.J.: Department of Mathematics,  University of
Illinois, Urbana,
 IL 61801, USA\\
 junge@math.uiuc.edu\\

\noindent Q.X.: Laboratoire de Math{\'e}matiques, Universit{\'e} de
France-Comt{\'e},
25030 Besan\c{c}on Cedex,  France\\
qxu@univ-fcomte.fr\\}

\end{document}